\documentclass[%
  a4paper
]{mpi2015-cscpreprint}

\usepackage{cleveref}

\usepackage{amsmath}
\usepackage{amssymb}
\usepackage{amsthm}

\usepackage[linesnumbered, lined, algoruled]{algorithm2e}

\crefname{algocf}{alg.}{algs.}
\Crefname{algocf}{Algorithm}{Algorithms}

\usepackage{subcaption}

\usepackage{pgfplots}
\pgfplotsset{compat=newest}
\usetikzlibrary{positioning}
\usetikzlibrary{patterns}
\usetikzlibrary{external}
\theoremstyle{definition}\newtheorem{remark}{Remark}
\theoremstyle{definition}\newtheorem{assumption}{Assumption}
\theoremstyle{plain}\newtheorem{proposition}{Proposition}
\theoremstyle{plain}\newtheorem{theorem}{Theorem}

\Crefname{assumption}{Assumption}{Assumptions}

\newcommand{\Hinf}{\ensuremath{\mathcal{H}_{\infty}}}
\newcommand{\LQG}{\ensuremath{{\mathrm{LQG}}}}
\newcommand{\trans}{\ensuremath{\mkern-1.5mu\mathsf{T}}}

\newcommand\RE{\mathsf{Re}}
\newcommand\infB{\mathcal{B}}
\newcommand\infC{\mathcal{C}}
\newcommand\vinf{\ensuremath{v_\infty}}
\newcommand\vd{\ensuremath{v_\delta}}
\def\yd{\ensuremath{y_\delta}}

\def\obsak{\widetilde{\widehat{A}}}
\def\obsbk{\widetilde{\widehat{B}}}
\def\obsck{\widetilde{\widehat{C}}}
\def\obsxk{\widetilde{\widehat{x}}}
\def\dobsxk{\dot{\widetilde{\widehat{x}}}}

\def\XPi{\ensuremath{X_\Pi}}
\def\YPi{\ensuremath{Y_\Pi}}

\def\Ainf{\ensuremath{A^{(\infty)}}}
\def\Aell{\ensuremath{A^{(\ell)}}}
\def\vell{\ensuremath{v^{(\ell)}}}
\def\Astokes{\ensuremath{A_{\mathsf{S}}}}

\def\dblcyl{\texttt{double cylinder}}
\def\cylwake{\texttt{cylinder wake}}

\newcommand{\R}{\ensuremath{\mathbb{R}}}
\newcommand{\C}{\ensuremath{\mathbb{C}}}

\newcommand{\nv}{\ensuremath{n_{\mathrm{v}}}}
\newcommand{\np}{\ensuremath{n_{\mathrm{p}}}}

\newcommand{\hA}{\ensuremath{\widehat{A}}}
\newcommand{\hB}{\ensuremath{\widehat{B}}}
\newcommand{\hC}{\ensuremath{\widehat{C}}}
\newcommand{\hE}{\ensuremath{\widehat{E}}}
\newcommand{\hG}{\ensuremath{\widehat{G}}}
\newcommand{\hK}{\ensuremath{\widehat{K}}}
\newcommand{\hM}{\ensuremath{\widehat{M}}}
\newcommand{\hN}{\ensuremath{\widehat{N}}}
\newcommand{\hX}{\ensuremath{\widehat{X}}}
\newcommand{\hY}{\ensuremath{\widehat{Y}}}
\newcommand{\hZ}{\ensuremath{\widehat{Z}}}
\newcommand{\hepsilon}{\ensuremath{\hat{\epsilon}}}
\newcommand{\hgamma}{\ensuremath{\hat{\gamma}}}

\newcommand{\tA}{\ensuremath{\widetilde{A}}}
\newcommand{\tB}{\ensuremath{\widetilde{B}}}
\newcommand{\tC}{\ensuremath{\widetilde{C}}}
\newcommand{\tE}{\ensuremath{\widetilde{E}}}
\newcommand{\tx}{\ensuremath{\tilde{x}}}

\newcommand{\cA}{\ensuremath{\mathcal{A}}}
\newcommand{\cB}{\ensuremath{\mathcal{B}}}
\newcommand{\cC}{\ensuremath{\mathcal{C}}}
\newcommand{\cE}{\ensuremath{\mathcal{E}}}
\newcommand{\cF}{\ensuremath{\mathcal{F}}}
\newcommand{\cL}{\ensuremath{\mathcal{L}}}

\newcommand{\opt}{\ensuremath{\mathsf{opt}}}
\def\gai{\ensuremath{\Gamma_{\text{in}}}}
\def\gao{\ensuremath{\Gamma_{\text{out}}}}
\def\ttol{\texttt{tol}}

\DeclareMathOperator{\diag}{diag}

\def\BT{balanced truncation}
\def\FEs{finite elements}

\definecolor{mpiblue}{HTML}{33a5c3}
\definecolor{mpigreen}{HTML}{007675}
\definecolor{mpired}{HTML}{78004B}
\definecolor{mpisand}{HTML}{ece9d4}

\begin{document}


\title{Robust output-feedback stabilization for incompressible flows using 
  low-dimensional \texorpdfstring{$\Hinf$}{H-infinity}-controllers}

\author[$\ast$]{Peter Benner}
\affil[$\ast$]{Max Planck Institute for Dynamics of Complex Technical
  Systems, Sandtorstra{\ss}e 1, 39106 Magdeburg, Germany.
  \email{benner@mpi-magdeburg.mpg.de}, \orcid{0000-0003-3362-4103}
  \authorcr \itshape
  Otto von Guericke University Magdeburg, Faculty of Mathematics,
  Universit{\"a}tsplatz 2, 39106 Magdeburg, Germany.
  \email{peter.benner@ovgu.de}}
  
\author[$\dagger$]{Jan Heiland}
\affil[$\dagger$]{Max Planck Institute for Dynamics of Complex Technical
  Systems, Sandtorstra{\ss}e 1, 39106 Magdeburg, Germany.
  \email{heiland@mpi-magdeburg.mpg.de}, \orcid{0000-0003-0228-8522}
  \authorcr \itshape
  Otto von Guericke University Magdeburg, Faculty of Mathematics,
  Universit{\"a}tsplatz 2, 39106 Magdeburg, Germany.
  \email{jan.heiland@ovgu.de}}
  
\author[$\ddag$]{Steffen W. R. Werner}
\affil[$\ddag$]{Max Planck Institute for Dynamics of Complex Technical
  Systems, Sandtorstra{\ss}e 1, 39106 Magdeburg, Germany.
  \email{werner@mpi-magdeburg.mpg.de}, \orcid{0000-0003-1667-4862}}

\shorttitle{Flow stabilization using low-order $\Hinf$-controllers}
\shortauthor{P. Benner, J. Heiland, S.~W.~R. Werner}
\shortdate{2021-09-29}
  
\keywords{robust control, incompressible flows, stabilizing feedback controller}

\msc{}
  
\abstract{%
  Output-based controllers are known to be fragile with respect to model
  uncertainties.
  The standard \Hinf-control theory provides a general approach to robust
  controller design based on the solution of the $\Hinf$-Riccati equations. 
  In view of stabilizing incompressible flows in simulations, two major
  challenges have to be addressed: the high-dimensional nature of the spatially
  discretized model and the differential-algebraic structure that comes with
  the incompressibility constraint. 
  This work demonstrates the synthesis of low-dimensional robust controllers
  with guaranteed robustness margins for the stabilization of incompressible
  flow problems.
  The performance and the robustness of the reduced-order controller with
  respect to linearization and model reduction errors are investigated and
  illustrated in numerical examples.
}

\maketitle


\section{Introduction}

We consider the incompressible Navier-Stokes equations with inputs and outputs
\begin{subequations}\label{eqn:intro-nse}
  \begin{align}
    \dot v &= -(v\cdot \nabla)v + \frac{1}{\RE} \Delta v - \nabla p + \infB u, \\
    0 &= \nabla \cdot v, \\
    y &= \infC v,
  \end{align}
\end{subequations}
and the question when a linear output-feedback controller $\mathcal{K}\colon y
\mapsto u$ can stabilize this nonlinear system around a possibly unstable steady
state in the presence of system uncertainties.
Here, the variables $v$ and $p$ describe the evolution of the velocity and
pressure fields in a given flow setup that is parametrized through the
Reynolds number $\RE$.
The operator $\infB$ models the actuation through the controls, and $\infC$ is
the output operator.

We will approach this question through a semi-discrete and linearized
approximation to~\cref{eqn:intro-nse}, model order reduction to cope with the
high dimensionality of the controller design problem, and the design of
controllers that can compensate for a large class of system uncertainties. 
Basically, our argument is that discretization and model reduction errors are of
the same nature such that a proven robustness margin can possibly overcome
unmodeled uncertainties, too.
Anyways, in order to potentially work in physical setups, any model-based
controller needs a certain robustness against inevitable model errors.
This rules out the standard \emph{linear quadratic Gaussian (LQG)} design that
has no guaranteed stability margin~\cite{Doy78}.
A general remedy is provided by $\Hinf$-controllers that, provably, can
compensate for linearization errors~\cite{BenH16,Hei21,BenHW19}, discretization
errors~\cite{BenH17,Cur06}, and truncation errors~\cite{morMusG91}.

The \Hinf-theory roots in the 1980s~\cite{Zam81}; see
also~\cite{FraD87,Fra87} for the historical background.
In view of its application in simulations, the development of state-space
formulations~\cite{DoyGKF89,FraD87} meant a breakthrough since it
came with general formulas for the controller design based on the solutions of
indefinite Riccati equations. 
Nonetheless, the computational effort for solving these Riccati equations is
significant so that, up to now, this design approach has rarely been considered
in large-scale simulations let alone the case where algebraic constraints are
present, that is for \emph{differential-algebraic equations} (DAEs) or
\emph{descriptor systems}.
If one leaves aside the \emph{offline} effort for designing the controller, the
theory seems well suited for large-scale problems since the controllers allow
for a low-dimensional approximation with a-priori estimates on performance and
robustness~\cite{morMusG91}.
The involved reduction is based on \emph{\BT} and the related \LQG-approach has
been investigated for descriptor systems in~\cite{morMoeRS11}.
Thus, with flow control in mind, the solution of large-scale Riccati
equations related to DAEs enables the use of the general \Hinf-theory.
It has been acknowledged that for general DAE systems, the standard
symmetric Riccati equations are only suitable under very restrictive
conditions~\cite{BenL87b}.
A nonsymmetric version has been shown to provide a true generalization for the
\Hinf-controller~\cite{BenH20} to the \emph{index-1} case and is widely
applicable for LQG-design for \emph{impulse controllable} descriptor 
systems~\cite{HeiZ21,WanYC98}.

In this work, we consider the incompressible Navier-Stokes equations with
control inputs in the momentum equation and velocity measurements only, so that
the input-to-output behavior can be equivalently realized as system of
\emph{ordinary differential equations (ODEs)}. 
Still, we keep the DAE structure since the corresponding transformation will
not be available in practice.
Then, the challenge is to realize what is suggested by the ODE theory without
explicitly resorting to transformations and projections.
Therefor, we adapt the established technique of realizing the projections
through the solution of saddle-point
systems~\cite{BaeBSetal15,Hei16,morHeiSS08,Wei16}, such that structure and
sparsity are preserved during all operations.

The $\Hinf$-feedback control for the 2D incompressible Navier-Stokes equations
has been treated in~\cite{DhaRT11} from a theoretical perspective.
Therein, unmodeled boundary inputs are considered and existence of optimal
feedback solutions is shown with the help of Riccati equations.
The general \Hinf-control problem is discussed in~\cite[Ch. 5]{Bar11a}.
Here, summing up the relevant research, the established state-space theory is
adapted to the infinite-dimensional incompressible Navier-Stokes equations. 

In this paper, we explore the Riccati-based \Hinf-controller design for
spatially discretized two-dimensional incompressible flows.
To this end,
\begin{itemize}
  \item we adapt the theory for the implicit treatment of the incompressibility
    constraint to the \Hinf-optimization problem,
  \item we leverage \Hinf-\BT{} to reduce the dimension of the controller design
    problem, 
  \item we provide numerically accessible formulas for a-priori estimation of
    the robustness of the controller with respect to the \Hinf-\BT{} model
    reduction error as well as linearization errors, and
  \item illustrate the performance in two challenging numerical examples.
\end{itemize}
As a result, we provide a complete numerical approach that makes
\Hinf-controller design feasible for large-scale Navier-Stokes systems and 
provides computable bounds on the robustness of the performance with respect to
both linearization errors~\cite{BenHW19} and model reduction
errors~\cite{morMusG91}.
The situation that the controller is based on inexact linearizations is relevant
in applications and fits well into the presented framework;
cf.~\cite{BenH17,BenHW19,Hei21}.
The presented numerical studies are based on the popular setup of the
two-dimensional wake behind a cylinder; 
see, for examples,~\cite{BerC08,GunL96,HeGMetal00,NoaAMetal03,Wil96}. 

This paper is organized as follows:
In \Cref{sec:basics}, we introduce the concepts of $\Hinf$-controller design and
truncation via the solution of Riccati equations and how the theory extends to
semi-discretized incompressible Navier-Stokes equations.
In \Cref{sec:methods}, we discuss numerical methods for the solution of
large-scale Riccati equations that implicitly  respect the incompressibility
constraint and provide a summary of steps for the design of robust
low-dimensional controllers with the accompanying formulas.
By means of two flow setups, we report on the performance of the resulting
low-dimensional controllers in \Cref{sec:examples}.
The paper is concluded in \Cref{sec:conclusions}.


\section{Mathematical basics}%
\label{sec:basics}
In this section, we introduce the basic concept and state-space approach of
robust \Hinf-controller design and how the relevant formulas are realized for
incompressible flow control setups.


\subsection{Riccati-based \texorpdfstring{$\Hinf$}{H-infinity} controller design}
\label{subsec:hinfcon}

For practical applications, the design of controllers that provide a certain
robustness against disturbances is essential, as neither models nor numerical
computations are fully able to match reality.
The $\Hinf$-control theory provides the design of such controllers.
In general, linear time-invariant systems of the form
\begin{subequations} \label{eqn:gensys}
\begin{align} 
  E \dot{x}(t) & = A x(t) + B_{1} w(t) + B_{2} u(t),\\
  z(t) & = C_{1} x(t) + D_{11} w(t) + D_{12} u(t),\\
  y(t) & = C_{2} x(t) + D_{21} w(t) + D_{22} u(t),
\end{align}
\end{subequations}
with $E, A \in \R^{n \times n}$, $B_{1} \in \R^{n \times m_{1}}$,
$B_{2} \in \R^{n \times m_{2}}$, $C_{1} \in \R^{p_{1} \times n}$,
$C_{2} \in \R^{p_{2} \times n}$, $D_{11} \in \R^{p_{1} \times m_{1}}$,
$D_{12} \in \R^{p_{1} \times m_{2}}$, $D_{21} \in \R^{p_{2} \times m_{1}}$,
and $D_{22} \in \R^{p_{2} \times m_{2}}$, are considered.
In~\cref{eqn:gensys}, the internal states $x$ are influenced by the control
inputs $u$ and disturbances $w$, and the user can observe the measurements
$y$ and performance outputs $z$ of the system.
The basic aim is to find a feedback controller $K\colon y \mapsto u$ that
internally stabilizes~\cref{eqn:gensys}.
Considering~\cref{eqn:gensys} in frequency domain allows the formulation of
the system's input-to-output relation in terms of its (partitioned with respect
to the different inputs and outputs) transfer function
\begin{align*}
  G(s) & = \begin{bmatrix} C_{1} \\ C_{2} \end{bmatrix} (sE - A)^{-1}
    \begin{bmatrix} B_{1} & B_{2} \end{bmatrix} +
    \begin{bmatrix} D_{11} & D_{12} \\ D_{21} & D_{22} \end{bmatrix}
    =: \begin{bmatrix} G_{11}(s) & G_{12}(s) \\ G_{21}(s) & G_{22}(s)
    \end{bmatrix}.
\end{align*}
Let also $K(s)$ denote the transfer function of the controller $K$, the
disturbance-to-performance behavior of the system can be formulated as
\begin{align} \label{eqn:lfractrafo}
  Z(s) & = (G_{11}(s) + G_{12}(s) K(s) (I_{p_{2}} - G_{22}(s) K(s))^{-1}
    G_{21}(s)) W(s) =: \cF(G, K) W(s),
\end{align}
where $Z$, $W$ are the Laplace transforms of the performances and disturbances,
and $I_{p_{2}}$ denotes the $p_{2}$-dimensional identity matrix.
$\cF$ is a \emph{lower linear fractional transformation}.
With~\cref{eqn:lfractrafo}, the \emph{optimal $\Hinf$-control problem} is
to find a controller $K$ that solves
\begin{align} \label{eqn:opthinf}
  \min\limits_{K~\text{stabilizing}} \lVert \cF(G, K) \rVert_{\Hinf} =:
    \gamma_{\opt},
\end{align}
where $\lVert.\rVert_{\Hinf}$ is the $\Hinf$-norm.
In general, this optimization problem~\cref{eqn:opthinf} is too difficult to
solve.
Instead, one can consider a relaxation in terms of the
\emph{sub-optimal $\Hinf$-control problem}: Find a stabilizing controller $K$
such that
\begin{align} \label{eqn:sopthinf}
  \lVert \cF(G, K) \rVert_{\Hinf} < \gamma,
\end{align}
which is then solved successively for decreasing robustness margins
$\gamma \rightarrow \gamma_{\opt}$.

A state-space theory has been developed that provides formulas for such
suboptimal controllers.
For ease of notation and theoretical derivations, we use the following set of
 assumptions:

\begin{assumption}[Simplifying and necessary assumptions for
  $\Hinf$-controller design]%
  \label{ass:sys-structure}
   In the formulation of system~\cref{eqn:gensys}, we have that
\begin{enumerate}
  \item $E$ is invertible,
  \item $D_{11}, D_{22} = 0$,
  \item $(sE - A, B_{1})$ is stabilizable and $(sE - A, C_{1})$ is detectable,
  \item $(sE - A, B_{2})$ is stabilizable and $(sE - A, C_{2})$ is detectable,
  \item $D_{12}^{\trans} \begin{bmatrix} C_{1} & D_{12} \end{bmatrix} = 
    \begin{bmatrix} 0 & I \end{bmatrix}$, and 
    $\begin{bmatrix} B_{1} \\ D_{21} \end{bmatrix} D_{21}^{\trans} = 
    \begin{bmatrix} 0 \\ I \end{bmatrix}$.
\end{enumerate}
\end{assumption}

We will comment on the assumptions in \Cref{rmk:assumptions} after we have
introduced the procedure to design an \Hinf-robust controller. 

\begin{proposition}[Existence of $\Hinf$-controllers~\cite{DoyGKF89,ZhoDG96}]%
  \label{prp:hinfexist}
  Given a $\gamma > \gamma_{\opt}$, there exists an admissible controller $K$ if
  and only if there are unique, symmetric positive semi-definite stabilizing
  solutions $X_{\Hinf}$ and $Y_{\Hinf}$ to the regulator and filter
  $\Hinf$-Riccati equations
  \begin{subequations} \label{eqn:hinfric}
  \begin{align} \label{eqn:hinfrric}
    A^{\trans} X E + E^{\trans} X A - E^{\trans} X 
    (B_{2}B_{2}^{\trans}-\gamma^{-2}B_{1}B_{1}^{\trans} ) 
    X E + C_{1}^{\trans} C_{1} & = 0,\\
    A Y E^{\trans} + E Y A^{\trans} - E Y
    (C_{2}^{\trans} C_{2} - \gamma^{-2} C_{1}^{\trans} C_{1}
      )YE^{\trans} + B_{1} B_{1}^{\trans} & = 0,
  \end{align}
  \end{subequations}
  and, additionally,
  \begin{align} \label{eqn:robustcond}
    \gamma^{2} > \lambda_{\max}(Y_{\Hinf} E^{\trans} X_{\Hinf} E),
  \end{align}
  where $\lambda_{\max}(M)$ denotes the maximum eigenvalue of the matrix $M$.
\end{proposition}

In the case that \Cref{prp:hinfexist} holds, a stabilizing controller
solving~\cref{eqn:sopthinf}, known in the literature as the \emph{central} or
\emph{minimum entropy controller}, is given via
\begin{align*}
  K:\left\{
  \begin{aligned}
    \tE \dot{\tx}(t) & = \tA \tx(t) + \tB y(t),\\
    u(t) & = \tC \tx(t),
  \end{aligned}\right.
\end{align*}
where the system matrices can be computed as
\begin{subequations} \label{eqn:hinf-controller-mats}
\begin{align}
  \tE & = E,\\
  \tA & = A + E Y_{\Hinf}(\gamma^{-2} C_{1}^{\trans}C_{1} -
    C_{2}^{\trans}C_{2}) - B_{2}B_{2}^{\trans} X_{\Hinf} E Z_{\Hinf},\\
  \tB & = E Y_{\Hinf} C_{2}^{\trans}, \\
  \tC & = -B_{2}^{\trans} X_{\Hinf} E Z_{\Hinf},
\end{align}
\end{subequations}
with $Z_{\Hinf} = (I_{n} - \gamma^{-2}Y_{\Hinf} E^{\trans} X_{\Hinf} E)^{-1}$;
see~\cite[Sec. III.C]{DoyGKF89} for the general case and also~\cite[Eqns. (16)
and (17)]{morMusG91} for the normalized case described in the following section.

\begin{remark}[Comments on the simplifying and necessary assumptions]%
  \label{rmk:assumptions}
  With $E$ invertible, our \Cref{ass:sys-structure} is equivalent to the
  standard set of assumptions made in \Hinf-robust controller design; cf., e.g.,
  \cite[Sec. IIA]{DoyGKF89}. 

  The assumptions on stabilizability and detectability with respect to $B_2$ and
  $C_2$ are clearly necessary for the existence of a stabilizing output feedback
  controller in the case of no disturbances.
  Stabilizability and detectability with respect to $B_1$ and $C_1$ are assumed
  to simplify both the formulas and the derivations, and can be relaxed towards
  conditions that ensure that the Riccati equations~\cref{eqn:hinfric} are
  solvable in the $\mathcal{H}_2$ case, i.e., as $\gamma \to \infty$.
  A further relaxation would possibly rule out the proposed Riccati-based
  approach; cf.~\cite[Ch. 17.1]{ZhoDG96}. 

  \Cref{ass:sys-structure}~(2.) is basically made for simplicity of the
  formulas.
  With the transformation $K = (I + D_{22}\hat K)^{-1}\hat K$ that maps a
  controller for $D_{22}=0$ onto an equivalent controller for the case that
  $D_{22}\neq 0$, the assumption of $D_{22}$ poses no restriction at all.
  Similarly, the restriction to $D_{11}= 0$ can be lifted but complicates the
  formulas a lot.

  Finally, the necessary part of \Cref{ass:sys-structure}~(5.) is that these
  matrices have full rank.
  The explicit structure can then be achieved by input and output
  transformations.
  If, however, the full-rank conditions is not met, then the formulation may
  lead to a singular control problem; see, e.g.,~\cite[Ch. 17.1]{ZhoDG96}.
\end{remark}


\subsection{The normalized \texorpdfstring{$\Hinf$}{H-infinity}-problem
  and low-rank robust controllers}

The target application of this work is output-based feedback control of a
nonlinear system that is robust against system uncertainties stemming from
linearization and reduction errors. 
As a general model for the linearization uncertainty, we will assume that the
system matrix $A$ is subjected to an additive perturbation. 
In this case, in a feedback arrangement, disturbance inputs will be induced by
the measurements of the perturbed state and enter the system as a perturbation
of the control input.
Accordingly, in \cref{eqn:gensys}, one can consider $B_{1} = B_{2} =:B$ and
$C_{1} = C_{2} =: C$, which leads to the so-called \emph{normalized}
\Hinf-control problem; see~\cite{morMusG91}.

A robust \Hinf-controller based on the \emph{normalized} problem then applies
in our context as follows. 
First, its robustness margin $\gamma$ can be weighed up against the system error
that arises from the truncation of certain states of the controller. 
Since a full-order controller would be of the same size as the system,
such a truncation significantly supports the efficient evaluation of the
feedback law during a simulation.
Second, the robustness can cover the system error that comes from an inaccurate
linearization.
In this section, we discuss both application scenarios and provide relevant
a-priori and a-posteriori estimates.

As for controller truncation, we consider the $\Hinf$-\BT{} approach as it has
been reported in~\cite{morMusG91}. 
Therefor, we introduce another special case of~\cref{eqn:gensys}, namely the
so-called \emph{normalized LQG system}
\begin{subequations} \label{eqn:normlqg}
\begin{align}
  E \dot{x}(t) & = Ax(t) + B w_{1}(t) + B u(t),\\
  z_{1}(t) & = C x(t),\\
  z_{2}(t) & = u(t),\\
  y(t) & = C x(t) + w_{2}(t).
\end{align}
\end{subequations}
In the terms of the general system~\cref{eqn:gensys}, the particular structural
assumptions made for~\cref{eqn:normlqg} mean that 
\begin{equation*}
  B_2 = B, \quad C_2=C, \quad B_1 = 
  \begin{bmatrix}
    B & 0
  \end{bmatrix}, \quad
  C_1 = 
  \begin{bmatrix}
    C \\ 0
  \end{bmatrix}
\end{equation*}
and 
\begin{equation*}
  D_{11} = 0, \quad D_{22}=0, \quad D_{21} = 
  \begin{bmatrix}
    0 & I
  \end{bmatrix}, \quad
  D_{12} = 
  \begin{bmatrix}
    0 \\ I
  \end{bmatrix},
\end{equation*}
and imply that \Cref{ass:sys-structure} is fulfilled if $(sE-A, B)$ is
stabilizable and $(sE-A, C)$ is detectable.
Moreover, this form of $B_1$, $B_2$, $C_1$ and $C_2$ implies that
\begin{align*}
  \begin{aligned}
    B_{1} B_{1}^{\trans}= B_{2} B_{2}^{\trans} & = B B^{\trans} &
      \text{and} && C_{1}^{\trans} C_{1} = C_{2}^{\trans} C_{2} & = C^{\trans} C
  \end{aligned}
\end{align*}
and, from \Cref{prp:hinfexist}, that the existence of an admissible controller,
which solves the suboptimal \Hinf-control problem~\cref{eqn:sopthinf} for the
gain $\gamma$, is equivalent to the existence of symmetric positive 
semi-definite matrices $X_{\Hinf}$ and $Y_{\Hinf}$, which solve
\begin{subequations} \label{eqn:hinfricnorm}
\begin{align} \label{eqn:hinffricnorm}
  A Y_{\Hinf} E^{\trans} + E Y_{\Hinf} A^{\trans} - (1 - \gamma^{-2}) E Y_{\Hinf}
    C ^{\trans}C Y_{\Hinf} E^{\trans} + B B ^{\trans} & = 0,\\ 
  \label{eqn:hinfrricnorm}
  A^{\trans} X_{\Hinf} E + E^{\trans} X_{\Hinf} A - (1 - \gamma^{-2}) E^{\trans} 
    X_{\Hinf} B B ^{\trans} X_{\Hinf} E + C^{\trans} C & = 0,
\end{align}
\end{subequations}
such that the spectrum condition~\cref{eqn:robustcond} holds and the matrix
pencils $sE - (A - (1 - \gamma^{-2}) E Y_{\Hinf} C ^{\trans}C)$ and $sE - (A -
(1 - \gamma^{-2}) B B ^{\trans} X_{\Hinf} E)$ are stable.

By means of the two matrices $X_{\Hinf}$ and $Y_{\Hinf}$, the so-called
$\Hinf$-\BT{} of the system can be computed~\cite[Prop. 4.10]{morMusG91} that
enables the truncation of systems or controller states with an a-priori control 
on the approximation error.
A practical implementation of this $\Hinf$-\BT{} (HINFBT) for large-scale systems
is shown in \Cref{alg:hinfbt} that uses the \emph{square root balancing approach}
on approximating low-rank factorizations of $X_{\Hinf}$ and $Y_{\Hinf}$.
An implementation of~\Cref{alg:hinfbt} for the dense system case can be
found in~\cite{morBenW19b}.

An error bound for the approximation computed by \Cref{alg:hinfbt} is given in
terms of normalized coprime factorizations:

\begin{proposition}[$\Hinf$-\BT{} error bound~\cite{morMusG91}]%
  \label{prp:hinferr}
  Let $G(s) = M(s)^{-1} N(s)$ and $\hG(s) = \hM(s)^{-1} \hN(s)$ be normalized
  left coprime factorizations of the original system $G$ and the reduced-order
  model $\hG$ obtained by HINFBT with the robustness margin $\gamma >
  \gamma_{\opt}$, respectively.
  Then, a bound for the approximation error is given by
  \begin{align} \label{eqn:hinferr}
    \left\lVert \begin{bmatrix} \beta(N - \hN) &\; M - \hM \end{bmatrix}
      \right\rVert_{\Hinf} & \leq 2 \sum\limits_{k = r + 1}^{n}
      \frac{\beta\sigma_{k}}{\sqrt{1 + \beta^{2} \sigma_{k}^{2}}},
  \end{align}
  with $\beta = \sqrt{1 - \gamma^{-2}}$ and the characteristic $\Hinf$-values
  $\sigma_{k}$.
\end{proposition}

\begin{algorithm}[t]
  \SetAlgoHangIndent{1pt}
  \DontPrintSemicolon
  \caption{(Approximate) $\Hinf$-\BT~square root method.}%
  \label{alg:hinfbt}
  
  \KwIn{$A, B, C, E$ from~\eqref{eqn:normlqg}, robustness margin $\gamma > 0$.}
  \KwOut{Matrices of the reduced-order system $\hA, \hB, \hC, \hE$.}
  
  Compute low-rank approximations $Y_{\Hinf} \approx
    R R^{\trans}$ and $X_{\Hinf} \approx L L^{\trans}$ to the unique stabilizing
    solutions of~\cref{eqn:hinfricnorm}.\;
    
  Compute the singular value decomposition
    \begin{align*}
      L^{\trans} E R & = \begin{bmatrix} U_{1} & U_{2} \end{bmatrix}
        \begin{bmatrix} \Sigma_{1} & \\ & \Sigma_{2} \end{bmatrix}
        \begin{bmatrix} V_{1}^{\trans} \\ V_{2}^{\trans} \end{bmatrix},
    \end{align*}
    with $\Sigma_{1} = \diag(\sigma_{1}, \ldots, \sigma_{r})$ containing the
    $r$ largest characteristic $\Hinf$-values.\;
  
  Construct the truncation matrices
    \begin{align*}
      \begin{aligned}
        W & = L U_{1} \Sigma_{1}^{-\frac{1}{2}} & \text{and} &&
          T & = R V_{1} \Sigma_{1}^{-\frac{1}{2}}.
      \end{aligned}
    \end{align*}\vspace{-\baselineskip}\;

  Compute the reduced-order model by
    \begin{align*}
      \begin{aligned}
        \hE & = W^{\trans} E T = I_{r}, &
          \hA & = W^{\trans} A T, &
          \hB & = W^{\trans} B, &
          \hC & = C T.
      \end{aligned}
    \end{align*}\vspace{-\baselineskip}\;
\end{algorithm}

An important consequence of \Cref{prp:hinferr} is that the robustness margin
$\gamma$ can be put in context with the approximation error of the reduced-order 
model and, consequently, the final reduced-order controller.
Thereby, one can estimate the size of the reduced-order controller, which is
needed to still stabilize the original system.

In general, let $G(s) = M(s)^{-1} N(s)$ and $\hG(s) = \hM(s)^{-1} \hN(s)$ be as 
in \Cref{prp:hinferr} the normalized left coprime factorizations of the original 
system and the reduced-order approximation computed by HINFBT, and define
\begin{align} \label{eqn:sclcoprimeerr}
  \beta \hepsilon := \left\lVert \begin{bmatrix} \beta(N - \hN) & \; M - \hM
    \end{bmatrix} \right\rVert_{\Hinf}
\end{align}
to be the scaled coprime factor error, with $\beta = \sqrt{1 - \gamma^{-2}}$ and
$\gamma$ used in HINFBT.
Then, it has been shown in~\cite[Cor. 5.5]{morMusG91} that a sufficient
condition for the reduced-order central controller $\hK$ based on $\hG$ to
stabilize the original system $G$ is
\begin{align} \label{eqn:stabbound}
  \hepsilon (\beta + \hgamma) < 1,
\end{align}
where $\hgamma := \lVert \cF(\hG, \hK) \rVert_{\Hinf}$.
The condition given in~\cref{eqn:stabbound} arises to be not overwhelmingly
practical since it can be only evaluated after the computation of the
reduced-order model and the computation of the $\Hinf$-norm of a large-scale
error system is needed.
A more practical alternative to~\cref{eqn:stabbound} can be derived from the
a-priori estimate in \Cref{prp:hinferr}.
Consider the HINFBT error bound~\cref{eqn:hinferr} to define
\begin{align*}
  \beta \epsilon & := 2 \sum\limits_{k = r + 1}^{n}
    \frac{\beta\sigma_{k}}{\sqrt{1 + \beta^{2} \sigma_{k}^{2}}},
\end{align*}
or, more precisely,
\begin{align} \label{eqn:epsbound}
  \epsilon & := 2 \sum\limits_{k = r + 1}^{n}
    \frac{\sigma_{k}}{\sqrt{1 + \beta^{2} \sigma_{k}^{2}}}.
\end{align}
Note that $\hepsilon \leq \epsilon$ necessarily holds.
Since $\hG$ is constructed by HINFBT (\Cref{alg:hinfbt}) it also holds that 
$\hgamma \leq \gamma$ for $\gamma$ the robustness margin used in
\Cref{alg:hinfbt}.
Therefor, a sufficient a-priori condition for stabilization of the full-order 
system by the reduced-order central controller is
\begin{align} \label{eqn:stabbound2}
  \epsilon (\beta + \gamma) < 1.
\end{align}

\begin{remark}[Performance loss of reduced-order controller~\cite{morMusG91}]
  We need to note that while a reduced-order controller $\hK$ is capable to
  guarantee the stabilization of the full-order plant using conditions
  like~\cref{eqn:stabbound} or~\cref{eqn:stabbound2}, it negatively influences
  the performance gain of the closed-loop system.
  In fact, one can show that if the full-order central
  controller~\cref{eqn:hinf-controller-mats} satisfies $\lVert \cF(G, K)
  \rVert_{\Hinf} < \gamma$, then for the reduced-order
  controller $\hK$ based on a HINFBT reduced-order model (\Cref{alg:hinfbt})
  using the performance gain $\gamma$ it holds that 
  \begin{align} \label{eqn:robustbound}
    \lVert \cF(G, \hK) \rVert_{\Hinf} & \leq
      \hgamma + \frac{\hepsilon(1 + \hgamma)(1 + \beta +
      \hgamma)}{1 - \hepsilon(\beta + \hgamma)}
      < \gamma + \frac{\epsilon(1 + \gamma)(1 + \beta +
      \gamma)}{1 - \epsilon(\beta + \gamma)},
  \end{align}
  where $\hepsilon$ and $\epsilon$ describe the approximation error and its
  a-priori bound as in~\cref{eqn:sclcoprimeerr} and~\cref{eqn:epsbound},
  $\hgamma = \lVert \cF(\hG, \hK) \rVert_{\Hinf} \leq \gamma$ and $\beta =
  \sqrt{1 - \gamma^{-2}}$.
  In practice, the bounds in~\cref{eqn:robustbound} quickly converge to
  $\hgamma$ or $\gamma$, respectively, for increasing reduced order.
\end{remark}

Another purpose of the robustness margin $\gamma$ is to ensure stability of the
closed-loop system in the presence of system uncertainties.
Therefor, the amount of uncertainties that is guaranteed to be compensated is
given in the following proposition.

\begin{proposition}[Stabilization of disturbed
  {systems~\cite[Cor. 3.7]{McFG90}}]%
  \label{prp:stabdist}
  Given the system $G(s) = M(s)^{-1} N(s)$ of the form~\cref{eqn:normlqg} and a
  stabilizing controller $K$ that solves the sub-optimal $\Hinf$-controller
  problem~\cref{eqn:sopthinf}, i.e., $\lVert \cF(G, K) \rVert_{\Hinf} < \gamma$
  holds for a given robustness margin $\gamma$.
  Let $G_{\Delta}(s) = M_{\Delta}(s)^{-1} N_{\Delta}(s)$ be another system in
  normalized LQG form~\cref{eqn:normlqg}, then $K$ is guaranteed to also
  stabilize $G_{\Delta}$ if
  \begin{align} \label{eqn:gamma_vs_cpferror}
    \left\lVert \begin{bmatrix} N - N_{\Delta} &\; M - M_{\Delta} \end{bmatrix}
      \right\lVert_{\Hinf} < \gamma^{-1}.
  \end{align}
  holds.
\end{proposition}

In view of stabilization of incompressible Navier-Stokes equations by linear 
output-feedback controllers, the following considerations are relevant with 
respect to the estimate~\cref{eqn:gamma_vs_cpferror}.
It has been shown that an error in the linearization used for controller design,
smoothly transfers to a coprime factor perturbation in the transfer function;
see~\cite{BenH16,Hei21}. 
Accordingly, with increasing accuracy in the computation of the linearization,
the difference becomes arbitrarily small such that, eventually, a robust
controller based on a numerically computed linearization will be able to
stabilize the system.

As compact overview and to end this section, we are going to collect the
a-priori-type conditions for the construction of reduced-order stabilizing
controllers into the following theorem.
Note that this theorem allows to adaptively choose the size of the reduced-order
controller depending on the stabilization of the original full-order system
and the amount of disturbances that shall be compensated by the reduced-order
controller.
 
\begin{theorem}[Sufficient a-priori conditions for stabilization of
  disturbed systems]%
  \label{thm:romstabdist}
  
  Given a system $G = M(s)^{-1} N(s)$ in normalized LQG form~\cref{eqn:normlqg}
  and a reduced-order system $\hG$ computed by HINFBT (\Cref{alg:hinfbt}) with
  the robustness margin $\gamma$.
  The central controller $\hK$ based on $\hG$ is guaranteed to stabilize the
  full-order system $G$ if
  \begin{align*}
    \epsilon (\beta + \gamma) < 1,
  \end{align*}
  where $\beta = \sqrt{1 - \gamma^{-2}}$, and $\epsilon$ is computed from the
  truncated characteristic $\Hinf$-values~\cref{eqn:epsbound}.
  Additionally, the reduced-order central controller $\hK$ is guaranteed to
  stabilize all systems $G_{\Delta}(s) = M_{\Delta}(s)^{-1} N_{\Delta}(s)$
  for which
  \begin{align*}
    \left\lVert \begin{bmatrix} N - N_{\Delta} &\; M - M_{\Delta} \end{bmatrix}
      \right\lVert_{\Hinf} < \big(\gamma_{\mathrm{G}\widehat{\mathrm{K}}}
      \big)^{-1}
  \end{align*}
  holds, and where
  \begin{align*}
    \gamma_{\mathrm{G}\widehat{\mathrm{K}}} = \gamma +
      \frac{\epsilon(1 + \gamma)(1 + \beta +\gamma)}
      {1 - \epsilon(\beta + \gamma)}.
  \end{align*}
\end{theorem}


\subsection{Semi-discretization and linearization of Navier-Stokes equations}

Next, we briefly discuss how a linear finite dimensional ODE system can be
derived as a base for the controller design for the nonlinear incompressible
Navier-Stokes equations.
Details of the discretization and the modeling of boundary control will be
discussed together with the numerical examples in \Cref{sec:examples}.

A \emph{\FEs} discretization of the incompressible Navier-Stokes equations
\cref{eqn:intro-nse} leads to the semi-discrete system:
\begin{subequations}\label{eqn:nse-semidisc}
  \begin{align}
  E \dot v &= \Astokes v + F(v) + J^{\trans}p + Bu + f, \\
  0 &= Jv + g, \\
  y &= Cv,
\end{align}
\end{subequations}
where $E \in \R^{\nv \times \nv}$ is the mass matrix,
$\Astokes \in \R^{\nv \times \nv}$ is the discrete approximation of
$\frac{1}{\RE}\Delta$, $F\colon \R^{\nv} \to \R^{\nv}$ models the
convection, $J \in \R^{\np \times \nv}$ and $J^{\trans}$ represent the discrete
divergence and gradient, $B \in \R^{\nv \times m}$ and $C \in \R^{p \times \nv}$
are the discretized input and output operators, and where $f \in \R^{\nv}$ and
$g \in \R^{\np}$ are the inhomogeneities that arise from the inclusion of the
inflow boundary condition in strong form.
Let $\vinf$ be the steady state with $u = 0$ so that with $\vd = v-\vinf$ and
$\Ainf := \Astokes + (\partial_v F)(\vinf)$,
\begin{subequations}\label{eqn:nse-linearized}
  \begin{align}
  E \dot \vd &= \Ainf\vd + J^{\trans}p + Bu, \\
  0 &= J\vd ,\\
  \yd &= C\vd,
\end{align}
\end{subequations}
provides a linearization that can be used for regulating the deviation $\vd$
from the steady state and, thus, for designing a controller for
stabilizing $\vinf$. 
With the assumption that the chosen \FEs~scheme is \emph{LBB-stable}, one can
define a discrete realization of the \emph{Leray projector}
\begin{equation}\label{eqn:leray-disc}
  \Pi^{\trans} := I_{\nv} - E^{-1}J^{\trans}(JE^{-1}J^{\trans})^{-1}J,
\end{equation}
that maps $v(t)$ into the kernel of $J$ along the orthogonal complement (in the
inner product induced by the mass matrix $E$) of $J^{\trans}$. 
Making use of $\Pi$ and the identities $\Pi E = E \Pi ^{\trans}$ -- which holds
for symmetric $E$ -- and $\Pi^{\trans}\vd = \vd$, we can eliminate the discrete
pressure $p$ and the algebraic constraint $0 = J\vd$ and
rewrite~\cref{eqn:nse-linearized} as ODE system
\begin{subequations} \label{eqn:nse-linearized-ode}
  \begin{align}
    E \dot \vd &= \Pi\Ainf\Pi^{\trans}\vd + \Pi Bu, \\
    \yd &= C\Pi^{\trans}\vd.
  \end{align}
\end{subequations}
For such a realization in terms of an ODE system~\cref{eqn:nse-linearized-ode},
\Hinf-robust controllers can be defined via the solutions to the projected
$\Hinf$-Riccati equations
\begin{subequations} \label{eqn:prjctd-hinf-ric}
\begin{align}
  \Pi \big(\Ainf\big)^{\trans}\Pi^{\trans}XE + E^{\trans}X \Pi \Ainf
    \Pi^{\trans} - (1 - \gamma^{-2}) E^{\trans} X \Pi B B^{\trans}
    \Pi^{\trans} X E + \Pi C^{\trans} C\Pi^{\trans} & = 0,\\
  \Pi \Ainf \Pi^{\trans} Y E^{\trans} + E Y \Pi \big(\Ainf\big)^{\trans}
    \Pi^{\trans} - (1 - \gamma^{-2}) E Y\Pi C^{\trans}C\Pi^{\trans}YE^{\trans} +
    \Pi  BB^{\trans}\Pi^{\trans} & = 0.
\end{align}
\end{subequations}

In the context of nonlinear flow control, the problem of linearization errors 
and compensation by reduced-order controllers has been considered
in~\cite{BenHW19}.
There it is mentioned that linearization errors in~\cref{eqn:nse-linearized}
from~\cref{eqn:nse-semidisc} result in additive disturbances on the transfer
function, which can be handled as disturbances in the coprime factorization.
In fact, the coprime factorization can be explicitly written down.

\begin{proposition}[Normalized left coprime factorization for flow
  problems~\cite{BenHW19}]%
  \label{prp:normcoprime}
  
  Given the system~\cref{eqn:nse-linearized} such
  that~\cref{eqn:nse-linearized-ode} restricted to the image of $\Pi$ is
  detectable.
  The normalized left coprime factorization of~\cref{eqn:nse-linearized}
  is given by
  \begin{align*}
    \begin{bmatrix} N(s) & M(s) \end{bmatrix} & = \cC (s \cE - \cA)^{-1}
      \begin{bmatrix} \cB & -\cL \end{bmatrix} +
      \begin{bmatrix} 0 & I_{p} \end{bmatrix},
  \end{align*}
  where the system matrices are
  \begin{align*}
    \begin{aligned}
      \cE & = \begin{bmatrix} E & 0 \\ 0 & 0 \end{bmatrix}, &
        \cA & = \begin{bmatrix} \Ainf - (1 - \gamma^{-2}) E Y_{\Hinf}C^{\trans}C
        & J^{\trans} \\ J & 0\end{bmatrix},\\
      \cC & = \begin{bmatrix} C & 0 \end{bmatrix}, &
        \cB & = \begin{bmatrix} B \\ 0 \end{bmatrix},\quad
        \cL = \begin{bmatrix} (1 - \gamma^{-2}) E Y_{\Hinf} C^{\trans} \\ 0
        \end{bmatrix},
    \end{aligned}
  \end{align*}
  and with $Y_{\Hinf}$, the stabilizing solution of~\cref{eqn:prjctd-hinf-ric}.
\end{proposition}

The representation in \Cref{prp:normcoprime} allows for the explicit computation 
of disturbances on the coprime factorizations, e.g., for the evaluation of the 
necessary robustness as in \cref{eqn:gamma_vs_cpferror}, accessible to
numerical simulations.


\section{Numerical methods}%
\label{sec:methods}


\subsection{Projector-free realization for incompressible flows}

In principle, the projected Riccati equations~\cref{eqn:prjctd-hinf-ric} could
be treated by established Riccati equation solvers.
However, in practice, the resulting system matrix $\Pi\Ainf\Pi^{\trans}$ will
be large scale and dense, which makes further computations barely feasible in
terms of computation time and memory consumption.
Also, a systematic error can easily be introduced by inaccurate computations
with~\cref{eqn:leray-disc}.
Therefore, like in many similar applications,
e.g.,~\cite{BaeBSetal15,morGugSW13,Hei16}, we derive a suitable implicit
implementation of the projection as it was proposed initially
in~\cite{morHeiSS08}.

We note that only the projected parts $X_\Pi:= \Pi^{\trans}X\Pi$ and
$Y_\Pi = \Pi^{\trans}Y\Pi$ of the Riccati solutions $X$ and $Y$ contribute to the
controller; see~\cref{eqn:hinf-controller-mats} and~\cite{BenH17b}; where
$X_\Pi$ and $Y_\Pi$ solve the equations
\begin{subequations} \label{eqn:prjctd-hinf-ric-prjctd-sols}
\begin{align} \label{eqn:prjctd-hinf-ric-reg}
  \Pi \big(\Ainf\big)^{\trans}\Pi^{\trans} \XPi E^{\trans} + E\XPi 
    \Pi\Ainf\Pi^{\trans} -(1 - \gamma^{-2}) E\XPi B B^{\trans} \XPi E^{\trans} +
    \Pi C^{\trans} C \Pi^{\trans} & = 0,\\
  \label{eqn:prjctd-hinf-ric-fil}
  \Pi \Ainf \Pi^{\trans} \YPi E + E^{\trans}\YPi \Pi \big(\Ainf\big)^{\trans}
    \Pi^{\trans} - (1 - \gamma^{-2}) E^{\trans} \YPi C^{\trans} C \YPi E +
    \Pi B B^{\trans} \Pi^{\trans} & = 0.
\end{align}
\end{subequations}
These equations~\cref{eqn:prjctd-hinf-ric-prjctd-sols} are derived
from~\cref{eqn:prjctd-hinf-ric} by pre- and postmultiplication
with $\Pi$ and $\Pi^{\trans}$, respectively, and by means of the
identities $\Pi^2 = \Pi$ and $\Pi E=E\Pi^{\trans}$.
In computational methods for large-scale sparse Riccati 
equations~\cite{BenBKetal18,BenLP08,Kle68,Sim16,SimSM14,Wei16},
low-rank factors $Z_k$ of the solution are computed such that, e.g.,
$\YPi \approx Z_k Z_k^{\trans}$, by applying repeated solves of shifted linear
systems that for~\cref{eqn:prjctd-hinf-ric-fil} read like
\begin{align} \label{eqn:shiftlin}
  (\Pi A \Pi^{\trans} + p_i E) Z & =  W.
\end{align}
With the requirement that the right-hand side lies in the appropriate subspace,
i.e., $W = \Pi^{\trans} W$, \cref{eqn:shiftlin} can be equivalently formulated
as 
\begin{align*}
  \begin{bmatrix} A + p_i E & J^{\trans} \\ J & 0 \end{bmatrix}
    \begin{bmatrix} Z \\ Z_\perp \end{bmatrix} & =
    \begin{bmatrix} E^{\trans} W \\ 0 \end{bmatrix}.
\end{align*}
Here, the scalar $p_i \in \C$ is a shift parameter that usually occurs in the 
considered iterative low-rank solvers like Newton-ADI or Krylov subspace 
methods; see, e.g.,~\cite{BenBKetal18,BenLP08,Kle68,Sim16,SimSM14,Wei16}; and
$Z_\perp$ is an auxiliary variable.


\subsection{Computation of low-rank controllers}
\label{sec:comp-lr-factors-cntrllrs}

The practical construction of a reduced-order controller
for~\cref{eqn:nse-semidisc} follows, in principle, the different steps mentioned
in this paper so far with some additional numerical tricks.
For simplicity, we give a final summary of the performed steps in the following:

\emph{Step 1. Computation of a suitable robustness margin:}
So far, the margin $\gamma$ was assumed to be given, but in fact it can be
computed utilizing the necessary and sufficient conditions for the existence
of a stabilizing controller from \Cref{prp:hinfexist}.
In practice, we solve repeatedly~\cref{eqn:prjctd-hinf-ric-prjctd-sols} for
different instances of the robustness margin and use~\cref{eqn:robustcond} to
determine the next iterate.
Thereby, a sufficient $\gamma$ can be computed.

\emph{Step 2. Construction of the reduced-order model:}
Now, the full-order system~\cref{eqn:nse-semidisc} needs to be reduced by
\Cref{alg:hinfbt}.
Therefore, we use the final low-rank solution factors $Z_{k}^{\mathsf{r}}$ and
$Z_{k}^{\mathsf{f}}$ of the projected algebraic Riccati
equations~\cref{eqn:prjctd-hinf-ric-reg} and~\cref{eqn:prjctd-hinf-ric-fil},
respectively, from the previous computation of the robustness margin in 
Step~1 corresponding to the computed $\gamma$ such that
$\XPi \approx Z_{k}^{\mathsf{r}} \left( Z_{k}^{\mathsf{r}} \right)^{\trans}$ and
$\YPi \approx Z_{k}^{\mathsf{f}} \left( Z_{k}^{\mathsf{f}} \right)^{\trans}$.
\Cref{alg:hinfbt} is then used for the system matrices $E$, $\Ainf$, $B$,
$C$ by setting the low-rank factors of the Riccati equations to be
$R = Z_{k}^{\mathsf{f}}$ and $L = Z_{k}^{\mathsf{r}}$.
The rest follows exactly \Cref{alg:hinfbt}.
The order $r$ of the reduced-order model can be determined by using the criteria
in \Cref{thm:romstabdist} for the stabilization of the full-order system and
possible disturbances.
This gives the reduced-order system matrices $\hE = I_{r}$, $\hA$, $\hB$ and
$\hC$.

\emph{Step 3. Construction of the reduced-order controller:}
The system matrices of the reduced-order central controller $\hK$ can then be
computed adapting the formulas in~\cref{eqn:hinf-controller-mats}.
First, an approximation to the solutions of the $\Hinf$-Riccati
equations~\cref{eqn:hinfric} for the reduced-order system is directly given by
\begin{align*}
  \begin{aligned}
    \hY_{\Hinf} & = W^{\trans} E Z_{k}^{\mathsf{f}}
      (Z_{k}^{\mathsf{f}})^{\trans} E^{\trans} W && \text{and} &
      \hX_{\Hinf} & = T^{\trans} E^{\trans} Z_{k}^{\mathsf{r}}
      (Z_{k}^{\mathsf{r}})^{\trans} E T,
  \end{aligned}
\end{align*}
with the low-rank factors $Z_{k}^{\mathsf{f}}$ and $Z_{k}^{\mathsf{r}}$ from the
computation of the robustness margin, and $W, T$ the truncation matrices from 
the HINFBT in \Cref{alg:hinfbt}.
Then, the system matrices of the reduced-order controller are given by
\begin{subequations}\label{eqn:reduced-cntrl-controller}
\begin{align}
  \widetilde{\hE} & = I_{r},\\
  \widetilde{\hA} & = \hA - (1 - \gamma^{-2}) \hY_{\Hinf} \hC^{\trans} \hC
    - \hB \hB^{\trans} \hX_{\Hinf} \hZ_{\Hinf},\\
  \widetilde{\hB} & = \hY_{\Hinf} \hC^{\trans},\\
    \widetilde{\hC} & = - \hB^{\trans} \hX_{\Hinf} \hZ_{\Hinf},
\end{align}
\end{subequations}
where $\hZ_{\Hinf} = (I_{r} - \gamma^2 \hY_{\Hinf} \hX_{\Hinf})^{-1}$.


\section{Numerical examples}%
\label{sec:examples}


\subsection{Example setups}

\begin{figure}
  \centering
  	\pgfmathsetmacro{\domfac}{7.5}
	\pgfmathsetmacro{\domx}{1.*\domfac}
	\pgfmathsetmacro{\domy}{.4*\domfac}
	\pgfmathsetmacro{\cylrad}{.05*\domfac}
	\pgfmathsetmacro{\cylcx}{.3*\domfac}
	\pgfmathsetmacro{\cylcy}{.2*\domfac}
	\begin{tikzpicture}
		\draw (\cylcx, \cylcy) circle (\cylrad);
		\draw [white,thick,domain=30:60] plot ({\cylcx+\cylrad*cos(\x)}, {\cylcy+\cylrad*sin(\x)});
		\draw [white,thick,domain=300:330] plot ({\cylcx+\cylrad*cos(\x)}, {\cylcy+\cylrad*sin(\x)});
		\node at (\cylcx+\cylrad,{\cylcy+\cylrad*cos(30)}) [above] {${\Gamma_1}$};
		\node at (\cylcx+\cylrad,{\cylcy-\cylrad*cos(30)}) [below] {${\Gamma_2}$};
		\node at (\cylcx-\cylrad,\cylcy) [left] {$\Gamma_w$};
		\node at (\cylcx+\cylrad,\cylcy) [right] {$\Gamma_w$};
		\draw (0,0) -- node[below]{$\Gamma_w$}(\domx,0) ;
		\draw (\domx,\domy) -- node[above]{$\Gamma_w$}(0,\domy);
		\draw [dashed] (\domx,0) -- node[right]{$\gao$}(\domx,\domy);
		\draw [dashed] (0,\domy) -- node[right]{$\gai$}(0,0);
	\end{tikzpicture}
  \caption{Computational domain of the cylinder wake.}
  \label{fig:cyldom}
\end{figure}
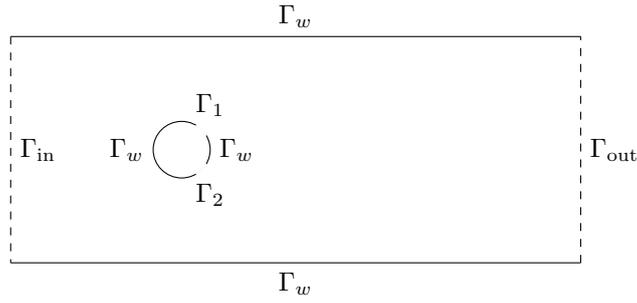

As numerical examples, we consider here two-dimensional flows through a channel 
with circular obstacles, exemplarily depicted in \Cref{fig:cyldom}, with
controls acting on the boundary of the obstacles, and with observation of 
locally spatially averaged velocities in a domain of observation downstream of 
the obstacles.

The generic model then reads
\begin{subequations}\label{eqn:numerics_generic_nse}
\begin{align}
  \dot v(t) +(v(t)\cdot \nabla)v(t) - \frac{1}{\RE} \Delta v(t) + \nabla
    p(t)&= 0, & \quad \text{in }&\Omega,\\
  \nabla \cdot v(t) &= 0, & \quad \text{in }&\Omega, \\
  \tfrac{1}{\RE}\tfrac{\partial v}{\partial n}v(t) -np(t) &= 0, & \quad
    \text{on } & \gao, \label{eqn:numgennse-outflow} \\
  v(t) &= 0, &  \quad \text{on } &\Gamma_w, \\
  v(t) &= ng_{\text{in}}, &  \quad \text{on } &\gai,\\
  v(t) &= ng_iu_i(t), &  \quad \text{on } &\Gamma_i, \quad i=1,2.
  \label{eqn:numggennse-controls}
\end{align}
\end{subequations}
Here, $n$ denotes the inward normal of the boundaries, $g_{\text{in}}$ models
the inflow boundary condition via a parabola that takes the value $0$ at the 
edges and $2/3$ at the center of the inflow boundary, respectively, and $g_i$ 
and $u_i$ are control shape functions and scalar control value functions.
The condition~\cref{eqn:numgennse-outflow} for the outflow is the standard
\emph{do-nothing} condition.
The geometrical extensions, the choice of the parameters
$u_{\text{in}}$ and $\RE$, and the definition of the output operator $\infC$,
and of the shape functions are given in the description of the test cases below
and in \Cref{tab:sim-setup}.


\subsubsection{Test case: \cylwake}

As first example, we consider the cylinder as it has been described in the
benchmark work~\cite{SchT96} in an enlarged domain to reduce the stabilizing
effects of the channel walls. 
To readily include the boundary controls in the \FEs~discretization, we relax
the Dirichlet control conditions~\cref{eqn:numggennse-controls} towards Neumann
conditions via
\begin{align} \label{eqn:nse-rob-cont-bcs}
  v(t) & = ng_i u_i(t) - \alpha (\tfrac{1}{\RE}\tfrac{\partial v}{\partial n}(t)
  -np(t)) \text{ on }\Gamma_i, \quad i=1,2,
\end{align}
with a parameter $\alpha$ that is supposed to be small; see, e.g.,~\cite{HouR98}
for convergence properties of this relaxation in optimal control of stationary
flows. 
The shape functions $g_i$ are chosen to be $g_i(x) = \cos(s_i(x)) - 1$, where
$s_i\colon \Gamma_i \to [0, 2\pi]$, for $i=1,2$, parametrizes the arc length of
the control boundaries; see~\cite[Section~9.3]{BehBH17}, where also the
assembling of the associated control operator is explained.

\begin{figure}[t]
  \centering
  \begin{subfigure}[b]{.495\textwidth}
    \centering
    \includegraphics[width=\textwidth]{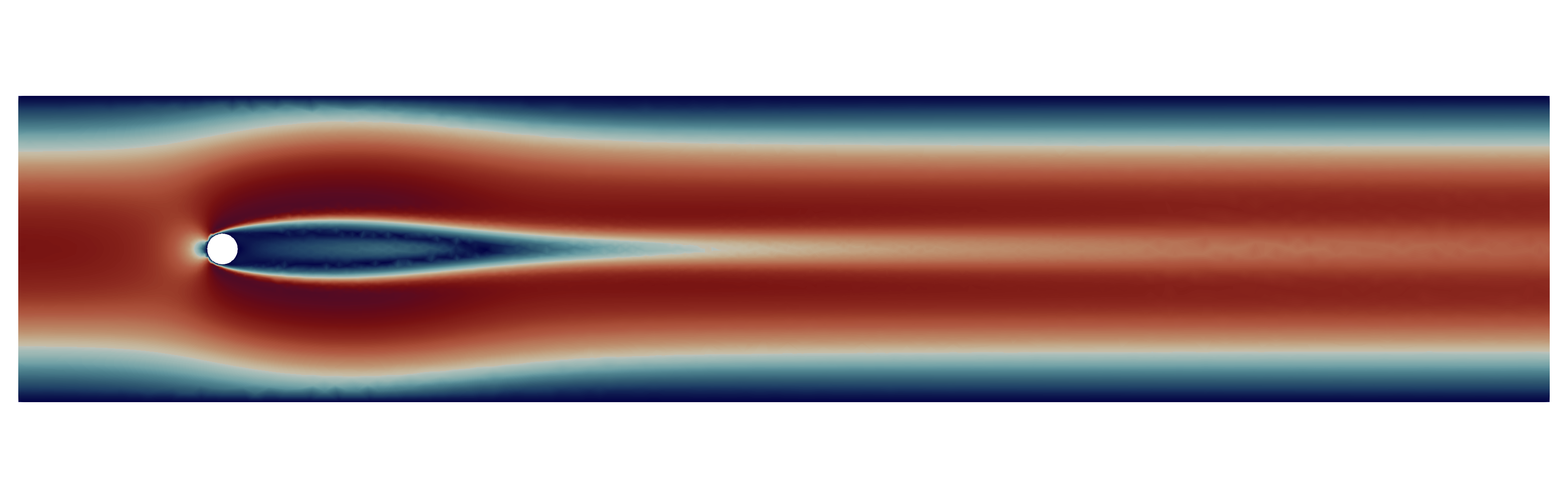}
    \vspace{-2\baselineskip}
    
    \caption{Steady state solution.}
  \end{subfigure}%
  \hfill%
  \begin{subfigure}[b]{.495\textwidth}
    \centering
    \includegraphics[width=\textwidth]{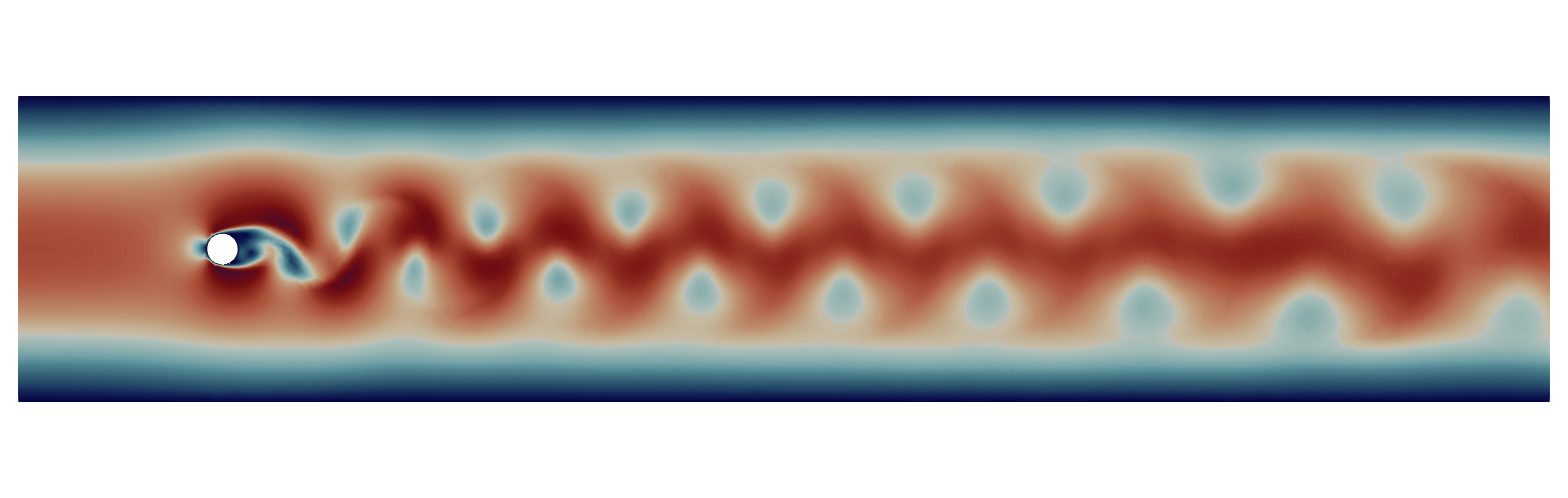}
    \vspace{-2\baselineskip}
    
    \caption{Snapshot of disturbed flow.}
  \end{subfigure}
  
  \begin{subfigure}{\textwidth}
    \centering
    \input{figures/cw60-nofb-lt.tikz}
  
    \caption{Measured output signal $y = C v$ over time $t$.
      The red lines depict the signal of the lateral averaged velocities,
      whereas the blue lines depict the measurements of the longitudinal
      velocities.}
  \end{subfigure}
  
  \caption{Example simulation of the cylinder wake with no control input at
    $\RE = 60$.}
  \label{fig:numsetup_cw}
\end{figure}


\subsubsection{Test case: \dblcyl}

\begin{figure}[t]
  \centering
  \begin{subfigure}[b]{.495\textwidth}
    \centering
    \includegraphics[width=\textwidth]{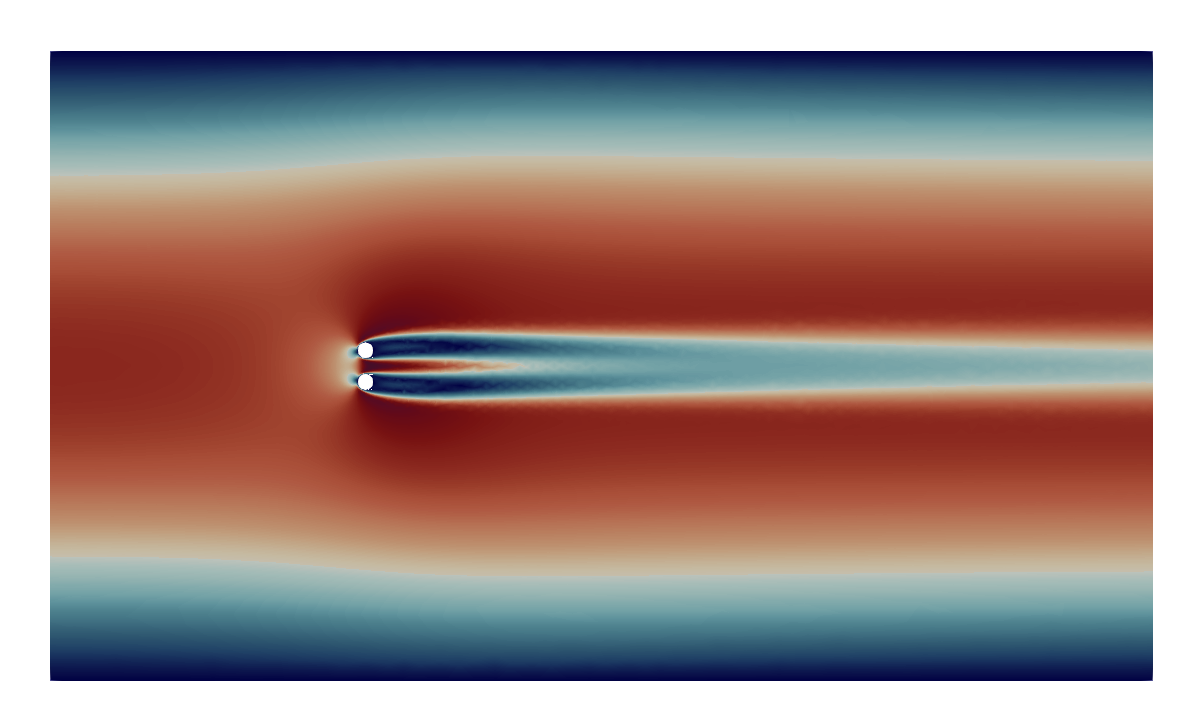}
    \vspace{-1.5\baselineskip}

    \caption{Steady state solution.}
  \end{subfigure}%
  \hfill%
  \begin{subfigure}[b]{.495\textwidth}
    \centering
    \includegraphics[width=\textwidth]{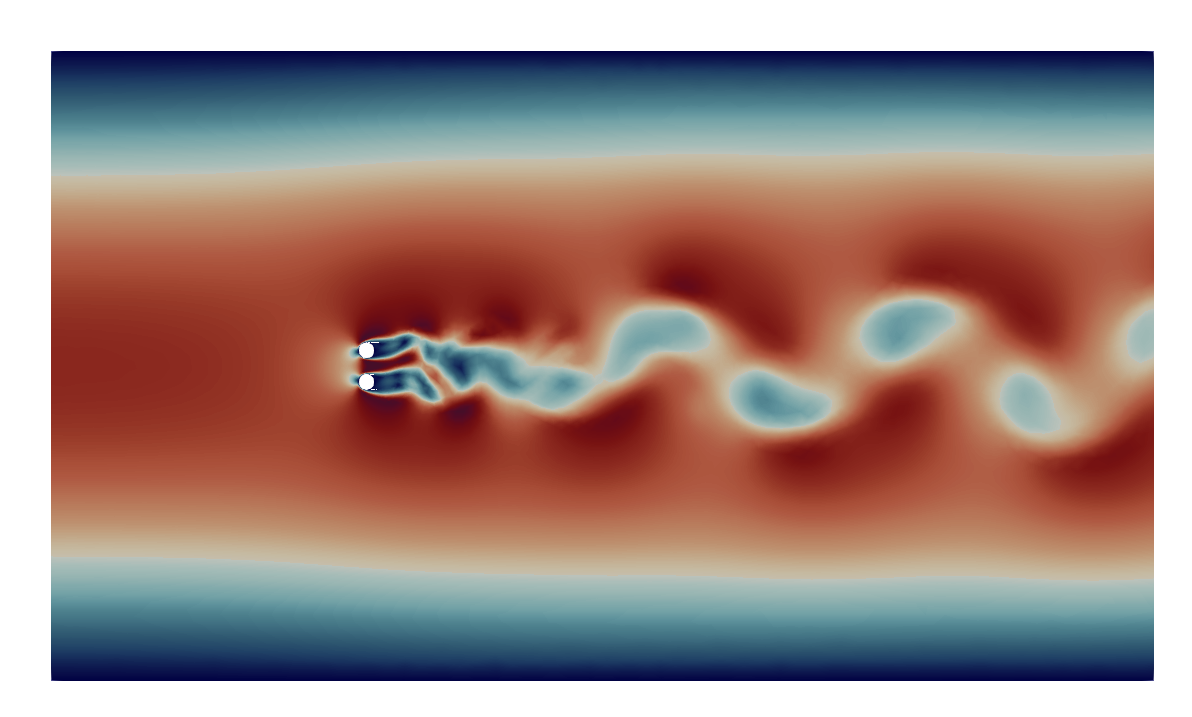}
    \vspace{-1.5\baselineskip}
    
    \caption{Snapshot of disturbed flow.}
  \end{subfigure}
  
  \begin{subfigure}{\textwidth}
    \centering
    \input{figures/drc60-nofb-lt.tikz}
  
    \caption{Measured output signal $y = C v$ over time $t$.
      The red lines depict the signal of the lateral averaged velocities,
      whereas the blue lines depict the measurements of the longitudinal
      velocities.}
  \end{subfigure}
  
  \caption{Example simulation of the double cylinder example with no control
    input at $\RE = 60$.}
  \label{fig:numsetup_dbrotcyl}
\end{figure}

As second example, we borrow the numerical setup from~\cite{morBorGZ16} of
the wake with two cylinders in two dimensions; see
\Cref{fig:numsetup_dbrotcyl}.
The actuation of the flow happens through controlled rotation of the individual
cylinders.
This means that, instead of the inflow
conditions~\cref{eqn:numggennse-controls}, we prescribe the control boundary
conditions as
\begin{align*}
  v\bigr|_{\Gamma_i} = t_{\Gamma_i}r_i u_i , \quad i=1,2\, ,
\end{align*}
where the $\Gamma_i$ are the two boundaries of the cylinders,
$t_{\Gamma_i}$ denote the tangential vectors at the boundaries, $r_i$ are the
radii, and the $u_i$ are the two scalar input functions, for $i=1,2$.
Although these boundary conditions have no normal component, e.g.,
$v\cdot n = 0$ at $\Gamma_i$, for $i=1,2$, and, thus, can be
included in a weak formulation with a bounded input operator $\infB$
(cf.~\cite[Sec. 3.1]{BarLT06}), we use the same Robin
relaxation~\cref{eqn:nse-rob-cont-bcs} that is defined to relax controls with
normal components.

\subsubsection{Further computational setup for both test cases}

The spatial discretization in both test cases is done with $P_2-P_1$
\emph{Taylor-Hood} \FEs.
The input operators are assembled as described in~\cite[Sec. 9.3]{BehBH17} with
the penalization parameter $\alpha$ set to $10^{-5}$. 
Since the output operator $\infC$ is of \emph{distributed type}, its assembling 
with \FEs~is straight forward.

We assemble the coefficients of the corresponding linearized
system~\cref{eqn:nse-linearized} and follow the procedure laid out in 
\Cref{sec:comp-lr-factors-cntrllrs}.
Therefor, we compute robustness margins $\gamma$ and the corresponding low-rank
approximations to the stabilizing solutions 
of~\cref{eqn:prjctd-hinf-ric-prjctd-sols} using the \emph{low-rank Riccati
iteration method} from~\cite{BenHW21}.
This approach is suited to solve Riccati equations like~\cref{eqn:hinfric}
with indefinite quadratic terms by splitting the computations into two steps:
first, the solution to the classical LQG-Riccati equations is computed and,
afterwards, residual Riccati equations are solved to update the overall
solution.
While this method was originally developed for Riccati equations with indefinite
quadratic terms, it also allows us here to efficiently compute the stabilizing
solutions of~\cref{eqn:prjctd-hinf-ric-prjctd-sols} for many different instances 
of the robustness margin to find suitable $\gamma$ for controller design.
Therein, we use the low-rank Newton-ADI~\cite{BenLP08,Wei16} as solver for the
LQG-Riccati equations and the low-rank RADI~\cite{BenBKetal18} for the update
residual equations.
The final computed robustness margins $\gamma$ for both test examples
are given in \Cref{tab:sim-params}.

The reduced central output-based feedback controller is then synthesized
via~\cref{eqn:reduced-cntrl-controller}, with its coefficients
\begin{align*}
  \begin{aligned}
    \obsak, && \obsbk, && \text{and}~\obsck.
  \end{aligned}
\end{align*}
Thus, with $\vinf$ denoting the target state and linearization point, the
closed-loop system reads
\begin{align*}
  E \dot v &= \Astokes v + F(v) + J^{\trans}p + B\obsck \obsxk +f, \\
  0 &= Jv + g, \\
  \dobsxk &= \obsak \obsxk + \obsbk C(v - \vinf).
\end{align*}

For the time discretization, we use a uniform grid of size $h$ and employ
\emph{backward differencing} of second order in the linear part including the
controller and the extrapolation 
\begin{align*}
  F(v(t+h)) & \approx 2F(v(t)) - F(v(t-h))
\end{align*}
for the nonlinearity.
Using one step of \emph{Heun's method} for the initialization, this results
in a time integration scheme of order two.
As the initial value, we use the corresponding steady state $v_\infty$.
To trigger the instabilities, we add an input perturbation that acts at the
beginning of the simulation like
\begin{equation*}
  u_{\delta} = 
  \begin{cases}
    \begin{aligned}
      10^{-6}\sin (2\pi t), & \quad \text{for}~t \in [0,1], \\
      0, & \quad \text{for}~t > 1.
    \end{aligned}
  \end{cases}
\end{equation*}
The final simulation setups and all parameters that define the simulations for
both test cases are listed in \Cref{tab:sim-setup,tab:sim-params}.

\begin{table}[t]
  \centering
  \caption{Simulation setups.
    For the computation of the Reynolds number, $\bar{v} = 1$ is the average
    inflow velocity and the dynamic viscosity $\nu$ is chosen accordingly.}
  \label{tab:sim-setup}
  \vspace{\baselineskip}
  
  {\def\arraystretch{1.5}%
  \begin{tabular}{c|p{5.2cm}|p{5.3cm}}
    & \cylwake & \dblcyl \\
    \hline
    Domain & $[0,5] \times [0,1]$ & $[-20,-20] \times [70,20]$\\
    Obstacles & Circle of radius $r_0=0.05$ \newline
      located at $(0.5, 0.67)$ &
      Two circles of radius $r_1, r_2 = 0.5$\newline
      located at $(0, 1)$ and $(0,-1)$\\
    Controls & Outlets at the cylinder periphery at $\pm \frac{\pi}3$ of
      arc length $\frac \pi 6$ &
      Independent rotation of both cylinders\\ 
    Observation domain & $[2,2.1] \times [0.3,0.7]$ horizontally split
      into $3$ equally sized subdomains &
      $[5, 6]\times[-2,2]$ horizontally split into $4$ equally sized
      subdomains\\
    Reynolds number $\RE$ & $60 = \frac{1\cdot 0.05}{\frac 56 10^{-3}}
      = \frac{\bar v \cdot r_0}{\nu}$ &
      $ 60 = \frac{1\cdot 1}{\frac 53 10^{-2}}
      = \frac{\bar v \cdot(r_1+r_2)}{\nu}$\\
  \end{tabular}}
\end{table}

\begin{table}
  \centering
  \caption{Simulation and controller parameters.}
  \label{tab:sim-params}
  \vspace{\baselineskip}

  {\def\arraystretch{1.5}%
  \begin{tabular}{c|p{3.5cm}|p{3.5cm}}
    & \cylwake & \dblcyl \\
    \hline
    Time interval & $[0,30]$& $[0, 300]$\\
    Time step size $h$ & $0.00075$ & $0.00390625$\\
    Relaxation parameter $\alpha$ & $10^{-5}$ &$10^{-5}$\\
    Dimension of $(v(t), p(t))$ & $(41718, 5418)$& $(45528, 5809)$\\
    Robustness margin $\gamma$ & $313.0176$ & $12.5418$\\
  \end{tabular}}
\end{table}

\begin{remark}[Validity of \Cref{ass:sys-structure}]
  The problem at hand derives from a state-space system with input operator
  $B$ and output operator $C$ that is brought into the normalized
  form~\cref{eqn:normlqg}. 
  Accordingly, \Cref{ass:sys-structure} is reduced to
  stability and detectability with respect to the given inputs and outputs; cf.
  the discussion after equation~\cref{eqn:normlqg}.
  Due to the involvement of the projector $\Pi$ in the underlying
  ODE~\cref{eqn:nse-linearized-ode}, for example, the system
  $(sE-\Pi \Ainf \Pi^{\trans}, \Pi B)$ will always have nonstabilizable modes on
  the imaginary axis so that it can not be stabilizable in the standard
  definition. 
  However, the modes associated with the kernel of $\Pi$ are excluded from the
  dynamics, so that notions of stabilizability and detectability can be adapted,
  e.g., by using the reduced coordinates as in~\cite{morHeiSS08}.

  Still, an analytical confirmation of \Cref{ass:sys-structure} is a difficult
  task.
  There exist relevant fundamental work (see, e.g.~\cite{NguR15}) but
  these results are generic and do not respect particular input operators and
  specific domains, let alone the discretization.
  Instead, one may resort to a numerical approach to establish stabilizability
  as proposed in~\cite[Sec. 5]{BenH16}.

  Finally, we want to remark that this stability and detectability is both
  necessary and sufficient for the existence of the stabilizing Riccati
  solutions and the convergence of the employed algorithms.
  Thus, a failure of the algorithms means that the assumptions do not hold and
  vice versa.
\end{remark}


\subsection{Numerical results}

The experiments reported here have been executed on a machine with 2 Intel(R)
Xeon(R) Silver 4110 CPU processors running at 2.10GHz and equipped with
192 GB total main memory.
The computer is run on CentOS Linux release 7.5.1804 (Core).

The solutions to the large-scale projected Riccati 
equations~\cref{eqn:prjctd-hinf-ric-prjctd-sols} have been computed in MATLAB 
9.7.0.1190202 (R2019b) using the routines from the M-M.E.S.S. library
version 2.0.1~\cite{morBenKS21, SaaKB20-mmess-2.0.1}.
Also, the errors in coprime factorizations have been computed in MATLAB.
For the spatial discretization, we employ the \emph{FEniCS}~\cite{fenics15}
\FEs~toolbox and the python module
\emph{dolfin\_navier\_scipy}~\cite{swHei-dns19} to extract the discrete
operators for the computation of the feedback gains and for the time
integration, which is done in \emph{SciPy}.

\begin{center}%
  \setlength{\fboxsep}{5pt}%
  \fbox{%
  \begin{minipage}{.92\textwidth}
    \textbf{Code and data availability}\newline
    The source code of the implementations used to compute the presented 
    results, the used data and the computed results are available from
    \begin{center}
      \href{https://doi.org/10.5281/zenodo.5532539}
        {\texttt{doi:10.5281/zenodo.5532539}}
    \end{center}
    under the MIT license and are authored by Jan Heiland and
    Steffen W. R. Werner.
  \end{minipage}}
\end{center}

We investigate the computed controllers in terms of robustness against
controller reduction and robustness against linearization errors via the
following criteria:
\begin{itemize}
  \item Does the robustness margin $\gamma$ cover the linearization error, i.e.,
    \begin{align*}
      \left\lVert \Delta \right\lVert_{\Hinf} :=
      \left\lVert \begin{bmatrix} N - N_{\Delta} &\; M - M_{\Delta}
        \end{bmatrix} \right\lVert_{\Hinf} & < \gamma^{-1}
    \end{align*}
    as in~\cref{eqn:gamma_vs_cpferror}?
  \item Does the robustness margin cover the a-priori estimate for the
    truncation error, i.e., $\epsilon( \beta + \gamma) < 1$ as
    in \Cref{thm:romstabdist}?
    (For reference, we also check the less conservative a-posteriori estimate
    $\hat{\epsilon}(\beta + \gamma) < 1$; cf.~\cref{eqn:stabbound}.)
  \item Does the controller work in the numerical experiment, i.e., does it
    stabilize the steady state by completely suppressing the oscillations that
    can be observed in the uncontrolled cases as illustrated in
    \Cref{fig:numsetup_cw,fig:numsetup_dbrotcyl}?
\end{itemize}

\begin{remark}[Comments on linearization errors]
  Our analysis of the linearization error follows the goal of robustifying the
  observer based controller.
  Thus, we will consider linearization errors like inexact computation of the
  linearization points as it may arise in numerical calculations. 
  Whether or not the linear controller is suitable for the
  nonlinear model is a question of performance, which is not considered here.
  In our numerical examples, we choose the linearization point as starting
  point and add small input perturbations such that we may well assume
  that we are close to the working point and the \emph{linearization principle
  for controller design} ensures the performance.
  For estimates on the admissible deviation from the working point in the
  infinite-dimensional Navier-Stokes model see, e.g.,~\cite{Ray06}.
\end{remark}

To examine the linearization error, we proceed as follows: 
For the \cylwake, we consider the linearization error that stems from an
incomplete iteration for the steady state computation. That is, in
\cref{eqn:nse-linearized-ode}, instead of the
\emph{exact} linearization $\Ainf=\Astokes + (\partial_v F)(\vinf)$ based on the
\emph{exact} steady state $\vinf$, we consider $\Aell := \Astokes + (\partial_v
F)(v_\ell)$, where $v_\ell$ is the approximation to $\vinf$ that is obtained
after $\ell$ \emph{Picard} iterations started from the steady state solution for
$\RE=40$.
For the \dblcyl, we found that at $\RE = 60$, the steady state is so unstable
that the \emph{Picard} iteration does not converge.
The computation can be done by a \emph{Newton} iteration that, however, does not
provide a smooth parametrization of the approximation in the relevant region
because of its fast convergence.
Therefor, we consider the perturbed coefficient $\Aell = \Astokes +
(\partial_v F)(\vell)$, where $\vell$ as the steady state solution to the
problem with the Reynolds number perturbed by $\ell$ thousandths, i.e.,
$\RE^{(\ell)} = (1 + \frac{\ell}{1000}) \cdot 60$. 
We parametrize the truncation error via the truncation threshold \ttol~that
defines $\ttol > \sigma_{r+1} \geq \sigma_{r+2} \geq ...$, i.e., the size of
those characteristic $\Hinf$-values that are discarded in \Cref{alg:hinfbt}. 

We have checked the performance in the simulations by examining the
empirical variances in the time series of the output signal in the third and
fourth quarters of the time interval.
If the difference between the fourth and third segment is negative, we
conclude that the oscillations were on the decline.
If the difference is positive but in the order of $10^{-15}$, we conclude that
the signals were dominated by numerical errors.
In both cases, the corresponding setup was reported as stabilizing in
the simulation.
In \Cref{tab:linerr}, we have tabulated the computed coprime factor errors of
the different linearizations and identified the threshold where
$\|\Delta_\ell\| < \gamma^{-1}$.

\begin{table}[t]
  \caption{Computed left coprime factor errors.}
  \label{tab:linerr}
  
  \begin{subtable}[t]{.45\textwidth}
    \centering
    \caption{Inexact linearizations for \cylwake{}.}
    \vspace{\baselineskip}
    
    \begin{tabular}{crr}
      $\ell$ & $\lVert \Delta^{(\ell)} \rVert_{\Hinf}$ & \\
      \hline\noalign{\medskip}
      $12$ & $1.1113$ & \\ \noalign{\medskip}
      $24$ & $0.0861$ & \\ \noalign{\medskip}
      $46$ & $0.0034$ & \\ \noalign{\medskip}
      $\boldsymbol{47}$ & $\boldsymbol{0.0029}$ & $\boldsymbol{< \gamma^{-1}}$
        \\ \noalign{\medskip}
      $48$ & $0.0029$ & \\ \noalign{\medskip}
      $96$ & $0.0007$ & \\
      \noalign{\medskip}\hline\noalign{\smallskip}
    \end{tabular}
  \end{subtable}%
  \hfill%
  \begin{subtable}[t]{.55\textwidth}
    \centering
    \caption{Perturbed $\RE$ numbers for \dblcyl{}.}
    \vspace{\baselineskip}
    
    \begin{tabular}{crr}
      $\ell$ & $\lVert \Delta^{(\ell)} \rVert_{\Hinf}$
        & \\
      \hline\noalign{\medskip}
      $-96$ & $0.3022$ & \\ \noalign{\medskip}
      $-48$ & $0.1540$ & \\ \noalign{\medskip}
      $-25$ & $0.0807$ & \\ \noalign{\medskip}
      $\boldsymbol{-24}$ & $\boldsymbol{0.0775}$ & $\boldsymbol{< \gamma^{-1}}$
        \\ \noalign{\medskip}
      $-12$ & $0.0389$ & \\ \noalign{\medskip}
      $-6$  & $0.0194$ & \\
      \noalign{\medskip}\hline\noalign{\smallskip}
    \end{tabular}
  \end{subtable}
\end{table}

The search of the parameter ranges for successful controller setups is
illustrated in \Cref{fig:perform-reg}.
Therein, each mark denotes a simulation that either failed (unstable behavior)
or has been stabilized.
For a better practical understanding, we added additionally to the truncation
tolerance \ttol{} the orders of the reduced controllers used in the simulations.
The a-priori and a-posteriori stability criteria are added as horizontal lines
and the bound for the linearization error as vertical line.
The intersection of both hatched areas are the regions where the controller 
setups are predicted to be stable by both criteria.

For the \cylwake, the estimates precisely confine the range where the 
controllers are functional; see \Cref{fig:perform-reg-cw}.
The a-priori estimate~\cref{eqn:stabbound2} turns out to be away from the
observed failures by a factor larger than $8$, and the bound on the
linearization error by a factor of $2^{\frac{1}{2}}$.
Successful stabilization has been observed further outside of the region
predicted by the bound on the linearization error and also outside
the a-posteriori estimate defined by the truncation error.

The simulations for the \dblcyl{} are shown in \Cref{fig:perform-reg-dc}.
Here, the picture is even more compliant with the predictions as the
a-posteriori estimate for the region of stabilization perfectly separates the
failed from the stabilized simulation.
The a-priori bound for the truncation error suggests a \ttol~of $0.002$, which
is a significant underestimate of the stability region and therefor not shown
in \Cref{fig:perform-reg-dc}.
On the other hand, the bound for the linearization is again too
conservative for the actual region of stabilization, this time by a factor of
$4$.
Curiously, increasing slightly the order of the \Hinf-controller allows to
perform successful simulations with a perturbation of the $\RE$ number up
to about $14\%$.
As for the \cylwake, this shows that successful stabilization can be observed
outside the region predicted by the linearization error bound while the
error estimates for the truncation error do not leave a margin.

\begin{figure}[t]
  \centering
  
  \begin{subfigure}[t]{\textwidth}
    \centering
    \begin{tikzpicture}
  \begin{loglogaxis}[
    mark repeat     = {1},
    width           = 10cm,
    height          = 7cm,
    cycle list name = color list,
    ylabel = {Truncation tolerance \ttol},
    log basis x=2,
    log basis y=2,
    xlabel = {Number of Picard iterations $\ell$},
    xmax = 57,
    xtick       = {17, 24, 34, 48},
    xticklabels = {17, 24, 34, 48},
    ytick       = {0.0064, 0.00166, 0.0004, 0.00021, 0.0001}, 
    yticklabels = {0.0064, 0.00166, 0.0004, 0.00021, 0.0001}, 
    grid=major,
    legend style = {legend cell align=left, legend pos = north east,
    font =\footnotesize},
    line width = 1pt,
    clip = false
 ]
 
 \addplot[
   color = mpired,
   only marks,
   mark=*,
   mark options={scale=3, fill=mpisand}
  ] table [y = trc, x = l] {figures/cw60sstbl.dat};
  
  \addplot[
    color = mpired,
    only marks,
    mark=otimes,
    mark options={scale=3, fill=mpigreen}
  ] table [y=trc, x=l] {figures/cw60sstbl.dat};

  \addplot[
    color=mpiblue, 
    only marks, 
    mark=*,
    mark options={scale=3, fill=mpisand}
  ] table [y=trc, x=l] {figures/cw60nstbl.dat};
  
  \addplot[
    color=mpiblue, 
    only marks, 
    mark=otimes,
    mark options={scale=3, fill=mpisand}
  ] table [y=trc, x=l] {figures/cw60nstbl.dat};
  
  \addplot[
    color=mpiblue, 
    only marks, 
    mark=*, 
    mark options={scale=2, fill=mpisand}
  ] table [y=trc, x=l] {figures/cw60nntbl.dat};
  
  \path[pattern = north east lines, pattern color = mpired!40]
     (axis description cs:0,0) rectangle
     ({rel axis cs:1,0}|-{axis cs:0,0.00166});
  \path[pattern=north west lines, pattern color=mpigreen!40]
    ({axis cs:47,0}|-{rel axis cs:0,0}) rectangle (axis description cs:1, 1);
    
  \draw[line width=3pt, dashed, mpired!70] 
    ({rel axis cs:0,0}|-{axis cs:0,0.00166}) --
    ({rel axis cs:1,0}|-{axis cs:0,0.00166});
  \node[fill=white, text=mpired!70] 
    at ({rel axis cs:.4,0}|-{axis cs:0,0.0009})
    {\large{$\hat \epsilon( \beta+ \gamma)<1$}};
    
  \draw[line width=3pt, dashed, mpired!70]
    ({rel axis cs:0,0}|-{axis cs:0,0.00021}) --
    ({rel axis cs:1,0}|-{axis cs:0,0.00021});
  \node[fill=white, text=mpired!70] 
    at ({rel axis cs:.4,0}|-{axis cs:0,0.00014})
    {\large{$\epsilon(\beta+\gamma)<1$}};
  
  \draw[line width=3pt, dashed, mpigreen!70]
    ({axis cs:47,0}|-{rel axis cs:0,0}) --
    ({axis cs:47,0}|-{rel axis cs:0,1});
  \node[fill=white, text=mpigreen!70, rotate=90] 
    at (axis description cs:.95,.6)
    {\large{$\gamma\|\Delta_\ell\|_\infty < 1$}};
  
  \node at (57, 0.0064) [anchor = west] {$30$};
  \node at (57, 0.0032) [anchor = west] {$33$};
  \node at (57, 0.0016) [anchor = west] {$35$};
  \node at (57, 0.0004) [anchor = west] {$40$};
  \node at (57, 0.0001) [anchor = west] {$48$};
  \node[rotate = 90] at (axis description cs:1.1,.5)
    {Reduced controller order $r$};
  \end{loglogaxis}
\end{tikzpicture}
    
    \caption{Performance regions for \cylwake.}
    \label{fig:perform-reg-cw}
  \end{subfigure}
  \vspace{\baselineskip}
  
  \begin{subfigure}[t]{\textwidth}
    \centering
    \begin{tikzpicture} 
  \begin{loglogaxis}[
    mark repeat     = {1},
    width           = 10cm,
    height          = 7cm,
    cycle list name = color list,
    ylabel = {Truncation tolerance \ttol},
    log basis x=2,
    log basis y=2,
    xlabel = {Disturbance of $\RE$ number $\ell$},
    x dir  = reverse,
    xtick       = {192, 96, 48, 24, 12, 6},
    xticklabels = {-192, -96, -48, -24, -12, -6},
    ytick       = {0.512, 0.128, 0.032, 0.008},
    yticklabels = {0.512, 0.128, 0.032, 0.008},
    grid=major,
    legend style = {legend cell align=left, legend pos = north east,
    font =\footnotesize},
    line width = 1pt,
    clip = false
 ]

  \addplot[
    color=mpiblue,
    only marks,
    mark=*,
    mark options={scale=2, fill=mpisand}
  ] table [y=trc, x=l] {figures/db60sntbl.dat};

  \addplot[
    color=mpired,
    only marks,
    mark=*,
    mark options={scale=3, fill=mpisand}
  ] table [y=trc, x=l] {figures/db60sstbl.dat};
        
  \addplot[
    color=mpired,
    only marks,
    mark=otimes,
    mark options={scale=3, fill=mpisand}
  ] table [y=trc, x=l] {figures/db60sstbl.dat};
  
  \addplot[
    color=mpiblue, 
    only marks, 
    mark=*,
    mark options={scale=3, fill=mpisand}
  ] table [y=trc, x=l] {figures/db60nstbl.dat};
  
  \addplot[
    color=mpiblue, 
    only marks, 
    mark=otimes,
    mark options={scale=3, fill=mpisand}
 ] table [y=trc, x=l] {figures/db60nstbl.dat};
 
  \addplot[
    color=mpiblue, 
    only marks, 
    mark=*, 
    mark options={scale=2, fill=mpisand}
  ] table [y=trc, x=l] {figures/db60nntbl.dat};

  \path[pattern=north east lines, pattern color=mpired!40]
    (axis description cs:0,0) rectangle
    ({rel axis cs:1,0}|-{axis cs:0,0.053});
  \path[pattern=north west lines, pattern color=mpigreen!40]
    ({axis cs:24,0}|-{rel axis cs:0,0}) rectangle (axis description cs:1, 1);
    
  \draw[line width=3pt, dashed, mpired!70] 
    ({rel axis cs:0,0}|-{axis cs:0,0.053}) --
    ({rel axis cs:1,0}|-{axis cs:0,0.053});
  \node[fill=white, text=mpired!70] 
    at ({rel axis cs:.3,0}|-{axis cs:0,0.016})
    {\large{$\hat \epsilon(\beta+\gamma)<1$}};
    

  \draw[line width=3pt, dashed, mpigreen!70]
    ({axis cs:24,0}|-{rel axis cs:0,0}) --
    ({axis cs:24,0}|-{rel axis cs:0,1});
  \node[fill=white, text=mpigreen!70] 
    at ({rel axis cs:.76,0}|-{axis cs:0,0.28})
    {\large{$\gamma\|\Delta_\ell\|_\infty < 1$}};

  \node at (4.3, 0.512) [anchor = west] {$12$};
  \node at (4.3, 0.128) [anchor = west] {$15$};
  \node at (4.3, 0.064) [anchor = west] {$16$};
  \node at (4.3, 0.045) [anchor = west] {$17$};
  \node at (4.3, 0.032) [anchor = west] {$18$};
  \node at (4.3, 0.008) [anchor = west] {$29$};
  \node[rotate = 90] at (axis description cs:1.1,.5)
    {Reduced controller order $r$};
  \end{loglogaxis}
\end{tikzpicture}
    
    \caption{Performance regions for \dblcyl.}
    \label{fig:perform-reg-dc}
  \end{subfigure}
  \vspace{.5\baselineskip}
  
  \begin{tikzpicture}
  \begin{axis}[
    hide axis,
    scale only axis,
    width = 1cm,
    xmin  = 0,
    xmax  = 1,
    ymin  = 0,
    ymax  = 1,
    legend columns = -1,
    legend style = {
      at     = {(0,0)},
      anchor = center,
      /tikz/every even column/.append style = {column sep = .3cm},
      font   =\footnotesize}
  ]
 
  \addplot[
    only marks,
    color        = mpiblue, 
    mark         = *, 
    mark options = {scale=2, fill=white},
    line width   = 1pt
  ] coordinates {(0,0)};
  \addlegendentry{Unstable simulation};
  
  \addplot[
    only marks,
    color        = mpiblue, 
    mark         = otimes, 
    mark options = {scale=3},
    line width   = 1pt
  ] coordinates {(0,0)};
  \addlegendentry{Stabilized simulation};
  
  \addplot[
    only marks,
    color        = mpired, 
    mark         = otimes, 
    mark options = {scale=3},
    line width   = 1pt
  ] coordinates {(0,0)};
  \addlegendentry{Stabilized simulation within computed bounds};
 
  \end{axis}
\end{tikzpicture}

  \caption{Performance of reduced controllers for the test examples.}
  \label{fig:perform-reg}
\end{figure}


\section{Conclusions}%
\label{sec:conclusions}

As the numerical examples have illustrated, the use of the general
Riccati-based low-order \Hinf-controller design is well feasible for
incompressible flows in simulations. 
The provided estimates on the guaranteed robustness have been proven reliable
though, in some cases, conservative.
Together with the theoretical results of earlier works that \Hinf-controller can
compensate various model errors, the availability of efficient general purpose
numerical methods is key for the applicability of these model-based controllers 
in simulations and even experiments.


\addcontentsline{toc}{section}{Acknowledgment}
\section*{Acknowledgment}
All authors have been supported by the German Research Foundation (DFG)
Research Training Group 2297 ``MathCoRe'', Magdeburg.


\addcontentsline{toc}{section}{References}

\bibliographystyle{abbrvurl}
\bibliography{exportref}

\begin{thebibliography}{10}

\bibitem{fenics15}
M.~S. Aln{\ae}s, J.~Blechta, J.~Hake, A.~Johansson, B.~Kehlet, A.~Logg,
  C.~Richardson, J.~Ring, M.~E. Rognes, and G.~N. Wells.
\newblock The {FEniCS} project version~1.5.
\newblock {\em Archive of Numerical Software}, 3(100):9--23, 2015.
\newblock \href {https://doi.org/10.11588/ans.2015.100.20553}
  {\path{doi:10.11588/ans.2015.100.20553}}.

\bibitem{BaeBSetal15}
E.~B{\"a}nsch, P.~Benner, J.~Saak, and H.~K. Weichelt.
\newblock {R}iccati-based boundary feedback stabilization of incompressible
  {N}avier-{S}tokes flows.
\newblock {\em {SIAM} J. Sci. Comput.}, 37(2):A832--A858, 2015.
\newblock \href {https://doi.org/10.1137/140980016}
  {\path{doi:10.1137/140980016}}.

\bibitem{Bar11a}
V.~Barbu.
\newblock {\em Stabilization of {N}avier-{S}tokes Flows}.
\newblock Communications and Control Engineering. Springer, London, 2011.
\newblock \href {https://doi.org/10.1007/978-0-85729-043-4}
  {\path{doi:10.1007/978-0-85729-043-4}}.

\bibitem{BarLT06}
V.~Barbu, I.~Lasiecka, and R.~Triggiani.
\newblock {\em Tangential boundary stabilization of {N}avier-{S}tokes
  equations}, volume 181 of {\em Mem. Am. Math. Soc.}
\newblock American Mathematical Society, 2006.
\newblock \href {https://doi.org/10.1090/memo/0852}
  {\path{doi:10.1090/memo/0852}}.

\bibitem{BehBH17}
M.~Behr, P.~Benner, and J.~Heiland.
\newblock Example setups of {N}avier-{S}tokes equations with control and
  observation: {S}patial discretization and representation via linear-quadratic
  matrix coefficients.
\newblock e-print arXiv:1707.08711, arXiv, 2017.
\newblock cs.MS.
\newblock URL: \url{https://arxiv.org/abs/1707.08711}.

\bibitem{BenL87b}
D.~J. Bender and A.~J. Laub.
\newblock The linear-quadratic optimal regulator for descriptor systems.
\newblock {\em {IEEE} Trans. Autom. Control}, 32(8):672--688, 1987.
\newblock \href {https://doi.org/10.1109/TAC.1987.1104694}
  {\path{doi:10.1109/TAC.1987.1104694}}.

\bibitem{BenBKetal18}
P.~Benner, Z.~Bujanovi{\'c}, P.~K{\"u}rschner, and J.~Saak.
\newblock {RADI}: a low-rank {ADI}-type algorithm for large scale algebraic
  {R}iccati equations.
\newblock {\em Numer. Math.}, 138(2):301--330, 2018.
\newblock \href {https://doi.org/10.1007/s00211-017-0907-5}
  {\path{doi:10.1007/s00211-017-0907-5}}.

\bibitem{BenH16}
P.~Benner and J.~Heiland.
\newblock Robust stabilization of laminar flows in varying flow regimes.
\newblock {\em IFAC-PapersOnLine}, 49(8):31--36, 2016.
\newblock 2nd {IFAC} Workshop on Control of Systems Governed by Partial
  Differential Equations {CPDE} 2016.
\newblock \href {https://doi.org/10.1016/j.ifacol.2016.07.414}
  {\path{doi:10.1016/j.ifacol.2016.07.414}}.

\bibitem{BenH17}
P.~Benner and J.~Heiland.
\newblock Convergence of approximations to {R}iccati-based boundary-feedback
  stabilization of laminar flows.
\newblock {\em IFAC-PapersOnLine}, 50(1):12296--12300, 2017.
\newblock 20th IFAC World Congress.
\newblock \href {https://doi.org/10.1016/j.ifacol.2017.08.2476}
  {\path{doi:10.1016/j.ifacol.2017.08.2476}}.

\bibitem{BenH17b}
P.~Benner and J.~Heiland.
\newblock Nonlinear feedback stabilization of incompressible flows via updated
  {R}iccati-based gains.
\newblock In {\em 2017 IEEE 56th Annual Conference on Decision and Control
  (CDC)}, pages 1163--1168, 2017.
\newblock \href {https://doi.org/10.1109/CDC.2017.8263813}
  {\path{doi:10.1109/CDC.2017.8263813}}.

\bibitem{BenH20}
P.~Benner and J.~Heiland.
\newblock Equivalence of {R}iccati-based robust controller design for index-1
  descriptor systems and standard plants with feedthrough.
\newblock In {\em 2020 European Control Conference (ECC)}, pages 402--407,
  2020.
\newblock \href {https://doi.org/10.23919/ECC51009.2020.9143771}
  {\path{doi:10.23919/ECC51009.2020.9143771}}.

\bibitem{BenHW19}
P.~Benner, J.~Heiland, and S.~W.~R. Werner.
\newblock Robust controller versus numerical model uncertainties for
  stabilization of {N}avier-{S}tokes equations.
\newblock {\em IFAC-PapersOnLine}, 52(2):25--29, 2019.
\newblock 3rd {IFAC/IEEE CSS} Workshop on Control of Systems Governed by
  Partial Differential Equation {CPDE} 2019.
\newblock \href {https://doi.org/10.1016/j.ifacol.2019.08.005}
  {\path{doi:10.1016/j.ifacol.2019.08.005}}.

\bibitem{BenHW21}
P.~Benner, J.~Heiland, and S.~W.~R. Werner.
\newblock A low-rank solution method for riccati equations with indefinite
  quadratic terms.
\newblock e-print, arXiv, 2021.
\newblock in preparation.

\bibitem{morBenKS21}
P.~Benner, M.~K{\"o}hler, and J.~Saak.
\newblock Matrix equations, sparse solvers: {M-M.E.S.S.}-2.0.1---{P}hilosophy,
  features and application for (parametric) model order reduction.
\newblock In P.~Benner, T.~Breiten, H.~Fa{\ss}bender, M.~Hinze, T.~Stykel, and
  R.~Zimmermann, editors, {\em Model Reduction of Complex Dynamical Systems},
  volume 171 of {\em International Series of Numerical Mathematics}, pages
  369--392. Birkh{\"a}user, Cham, 2021.
\newblock \href {https://doi.org/10.1007/978-3-030-72983-7_18}
  {\path{doi:10.1007/978-3-030-72983-7_18}}.

\bibitem{BenLP08}
P.~Benner, J.-R. Li, and T.~Penzl.
\newblock Numerical solution of large-scale {L}yapunov equations, {R}iccati
  equations, and linear-quadratic optimal control problems.
\newblock {\em Numer. Lin. Alg. Appl.}, 15(9):755--777, 2008.
\newblock \href {https://doi.org/10.1002/nla.622} {\path{doi:10.1002/nla.622}}.

\bibitem{morBenW19b}
P.~Benner and S.~W.~R. Werner.
\newblock {MORLAB} -- {Model Order Reduction LABoratory} (version 5.0), Aug.
  2019.
\newblock see also: \url{https://www.mpi-magdeburg.mpg.de/projects/morlab}.
\newblock \href {https://doi.org/10.5281/zenodo.3332716}
  {\path{doi:10.5281/zenodo.3332716}}.

\bibitem{BerC08}
M.~Bergmann and L.~Cordier.
\newblock Optimal control of the cylinder wake in the laminar regime by
  trust-region methods and {POD} reduced-order models.
\newblock {\em J. Comput. Phys.}, 227(16):7813--7840, 2008.
\newblock \href {https://doi.org/10.1016/j.jcp.2008.04.034}
  {\path{doi:10.1016/j.jcp.2008.04.034}}.

\bibitem{morBorGZ16}
J.~Borggaard, S.~Gugercin, and L.~Zietsman.
\newblock Feedback stabilization of fluids using reduced-order models for
  control and compensator design.
\newblock In {\em 2016 IEEE 55th Conference on Decision and Control (CDC)},
  pages 7579--7585, 2016.
\newblock \href {https://doi.org/10.1109/CDC.2016.7799440}
  {\path{doi:10.1109/CDC.2016.7799440}}.

\bibitem{Cur06}
R.~F. Curtain.
\newblock A robust {LQG}-controller design for {DPS}.
\newblock {\em Internat. J. Control}, 79(2):162--170, 2006.
\newblock \href {https://doi.org/10.1080/00207170500512985}
  {\path{doi:10.1080/00207170500512985}}.

\bibitem{DhaRT11}
S.~Dharmatti, J.-P. Raymond, and L.~Thevenet.
\newblock ${H}^\infty$ feedback boundary stabilization of two-dimensional
  {N}avier-{S}tokes equations.
\newblock {\em {SIAM} J. Control Optim.}, 49(6):2318--2348, 2011.
\newblock \href {https://doi.org/10.1137/100782607}
  {\path{doi:10.1137/100782607}}.

\bibitem{Doy78}
J.~Doyle.
\newblock Guaranteed margins for {LQG} regulators.
\newblock {\em {IEEE} Trans. Autom. Control}, 23(4):756--757, 1978.
\newblock \href {https://doi.org/10.1109/TAC.1978.1101812}
  {\path{doi:10.1109/TAC.1978.1101812}}.

\bibitem{DoyGKF89}
J.~Doyle, K.~Glover, P.~P. Khargonekar, and B.~A. Francis.
\newblock State-space solutions to standard $\mathcal{H}_2$ and
  $\mathcal{H}_{\infty}$ control problems.
\newblock {\em {IEEE} Trans. Autom. Control}, 34(8):831--847, 1989.
\newblock \href {https://doi.org/10.1109/9.29425} {\path{doi:10.1109/9.29425}}.

\bibitem{Fra87}
B.~A. Francis.
\newblock {\em A Course in $\mathcal{H}_{\infty}$ Control Theory}, volume~88 of
  {\em Lect. Notes Control Inf. Sci.}
\newblock Springer-Verlag, Berlin, Heidelberg, 1987.
\newblock \href {https://doi.org/10.1007/BFb0007371}
  {\path{doi:10.1007/BFb0007371}}.

\bibitem{FraD87}
B.~A. Francis and J.~C. Doyle.
\newblock Linear control theory with an {$H_{\infty}$} optimality criterion.
\newblock {\em {SIAM} J. Control Optim.}, 25(4):815--844, 1987.
\newblock \href {https://doi.org/10.1137/0325046} {\path{doi:10.1137/0325046}}.

\bibitem{morGugSW13}
S.~Gugercin, T.~Stykel, and S.~Wyatt.
\newblock Model reduction of descriptor systems by interpolatory projection
  methods.
\newblock {\em {SIAM} J. Sci. Comput.}, 35(5):B1010--B1033, 2013.
\newblock \href {https://doi.org/10.1137/130906635}
  {\path{doi:10.1137/130906635}}.

\bibitem{GunL96}
M.~D. Gunzburger and H.~C. Lee.
\newblock Feedback control of {K}arman vortex shedding.
\newblock {\em J. Appl. Mech.}, 63(3):828--835, 1996.
\newblock \href {https://doi.org/10.1115/1.2823369}
  {\path{doi:10.1115/1.2823369}}.

\bibitem{HeGMetal00}
J.-W. He, R.~Glowinski, R.~Metcalfe, A.~Nordlander, and J.~Periaux.
\newblock Active control and drag optimization for flow past a circular
  cylinder.
\newblock {\em J. Comput. Phys.}, 163(1):83--117, 2000.
\newblock \href {https://doi.org/10.1006/jcph.2000.6556}
  {\path{doi:10.1006/jcph.2000.6556}}.

\bibitem{Hei16}
J.~Heiland.
\newblock A differential-algebraic {R}iccati equation for applications in flow
  control.
\newblock {\em {SIAM} J. Control Optim.}, 54(2):718--739, 2016.
\newblock \href {https://doi.org/10.1137/151004963}
  {\path{doi:10.1137/151004963}}.

\bibitem{swHei-dns19}
J.~Heiland.
\newblock dolfin\_navier\_scipy: a python {S}cipy {FE}ni{CS} interface.
\newblock Github/Zenodo, June 2019.
\newblock \url{https://github.com/highlando/dolfin_navier_scipy}.
\newblock \href {https://doi.org/10.5281/zenodo.3238622}
  {\path{doi:10.5281/zenodo.3238622}}.

\bibitem{Hei21}
J.~Heiland.
\newblock Convergence of coprime factor perturbations for robust stabilization
  of {O}seen systems.
\newblock {\em Math. Control Relat. Fields}, 2021.
\newblock \href {https://doi.org/10.3934/mcrf.2021043}
  {\path{doi:10.3934/mcrf.2021043}}.

\bibitem{HeiZ21}
J.~Heiland and E.~Zuazua.
\newblock Classical system theory revisited for turnpike in standard state
  space systems and impulse controllable descriptor systems.
\newblock Technical Report 2007.13621, arxiv, 2020.
\newblock math.OC, accepted for publication in SIAM J. Control Optim.
\newblock URL: \url{http://arxiv.org/abs/2007.13621}.

\bibitem{morHeiSS08}
M.~Heinkenschloss, D.~C. Sorensen, and K.~Sun.
\newblock Balanced truncation model reduction for a class of descriptor systems
  with application to the {O}seen equations.
\newblock {\em {SIAM} J. Sci. Comput.}, 30(2):1038--1063, 2008.
\newblock \href {https://doi.org/10.1137/070681910}
  {\path{doi:10.1137/070681910}}.

\bibitem{HouR98}
L.~S. Hou and S.~S. Ravindran.
\newblock A penalized {N}eumann control approach for solving an optimal
  {D}irichlet control problem for the {N}avier-{S}tokes equations.
\newblock {\em {SIAM} J. Control Optim.}, 36(5):1795--1814, 1998.
\newblock \href {https://doi.org/10.1137/S0363012996304870}
  {\path{doi:10.1137/S0363012996304870}}.

\bibitem{Kle68}
D.~L. Kleinman.
\newblock On an iterative technique for {R}iccati equation computations.
\newblock {\em {IEEE} Trans. Autom. Control}, 13(1):114--115, 1968.
\newblock \href {https://doi.org/10.1109/TAC.1968.1098829}
  {\path{doi:10.1109/TAC.1968.1098829}}.

\bibitem{McFG90}
D.~C. McFarlane and K.~Glover.
\newblock {\em Robust Controller Design Using Normalized Coprime Factor Plant
  Descriptions}, volume 138 of {\em Lect. Notes Control Inf. Sci.}
\newblock Springer, Berlin, Heidelberg, 1990.
\newblock \href {https://doi.org/10.1007/BFB0043199}
  {\path{doi:10.1007/BFB0043199}}.

\bibitem{morMoeRS11}
J.~M{\"o}ckel, T.~Reis, and T.~Stykel.
\newblock Linear-quadratic {G}aussian balancing for model reduction of
  differential-algebraic systems.
\newblock {\em Internat. J. Control}, 84(10):1627--1643, 2011.
\newblock \href {https://doi.org/10.1080/00207179.2011.622791}
  {\path{doi:10.1080/00207179.2011.622791}}.

\bibitem{morMusG91}
D.~Mustafa and K.~Glover.
\newblock Controller reduction by $\mathcal{H}_\infty$-balanced truncation.
\newblock {\em {IEEE} Trans. Autom. Control}, 36(6):668--682, 1991.
\newblock \href {https://doi.org/10.1109/9.86941} {\path{doi:10.1109/9.86941}}.

\bibitem{NguR15}
P.~A. Nguyen and J.-P. Raymond.
\newblock Boundary stabilization of the {N}avier--{S}tokes equations in the
  case of mixed boundary conditions.
\newblock {\em {SIAM} J. Control Optim.}, 53(5):3006--3039, 2015.
\newblock \href {https://doi.org/10.1137/13091364X}
  {\path{doi:10.1137/13091364X}}.

\bibitem{NoaAMetal03}
B.~R. Noack, K.~Afanasiev, M.~Morzy{\'n}ski, G.~Tadmor, and F.~Thiele.
\newblock A hierarchy of low-dimensional models for the transient and
  post-transient cylinder wake.
\newblock {\em J. Fluid Mech.}, 497:335--363, 2003.
\newblock \href {https://doi.org/10.1017/S0022112003006694}
  {\path{doi:10.1017/S0022112003006694}}.

\bibitem{Ray06}
J.-P. Raymond.
\newblock Feedback boundary stabilization of the two-dimensional
  {N}avier-{S}tokes equations.
\newblock {\em {SIAM} J. Control Optim.}, 45(3):790--828, 2006.
\newblock \href {https://doi.org/10.1137/050628726}
  {\path{doi:10.1137/050628726}}.

\bibitem{SaaKB20-mmess-2.0.1}
J.~Saak, M.~K{\"o}hler, and P.~Benner.
\newblock {M-M.E.S.S.-2.0.1} -- {T}he {M}atrix {E}quations {S}parse {S}olvers
  library, Feb. 2020.
\newblock see also:~\url{https://www.mpi-magdeburg.mpg.de/projects/mess}.
\newblock \href {https://doi.org/10.5281/zenodo.3606345}
  {\path{doi:10.5281/zenodo.3606345}}.

\bibitem{SchT96}
M.~Sch{\"a}fer, S.~Turek, F.~Durst, E.~Krause, and R.~Rannacher.
\newblock Benchmark computations of laminar flow around a cylinder.
\newblock In E.~H. Hirschel, editor, {\em Flow Simulation with High-Performance
  Computers {II}: {DFG} priority research program results 1993--1995},
  volume~52 of {\em Notes Numer. Fluid Mech.}, pages 547--566. Vieweg,
  Wiesbaden, 1996.
\newblock \href {https://doi.org/10.1007/978-3-322-89849-4_39}
  {\path{doi:10.1007/978-3-322-89849-4_39}}.

\bibitem{Sim16}
V.~Simoncini.
\newblock Analysis of the rational {K}rylov subspace projection method for
  large-scale algebraic {R}iccati equations.
\newblock {\em {SIAM} J. Matrix Anal. Appl.}, 37(4):1655--1674, 2016.
\newblock \href {https://doi.org/10.1137/16M1059382}
  {\path{doi:10.1137/16M1059382}}.

\bibitem{SimSM14}
V.~Simoncini, D.~B. Szyld, and M.~Monsalve.
\newblock On two numerical methods for the solution of large-scale algebraic
  {R}iccati equations.
\newblock {\em {IMA} J. Numer. Anal.}, 34(3):904--920, 2014.
\newblock \href {https://doi.org/10.1093/imanum/drt015}
  {\path{doi:10.1093/imanum/drt015}}.

\bibitem{WanYC98}
H.-S. Wang, C.-F. Yung, and F.-R. Chang.
\newblock Bounded real lemma and ${H}_\infty$ control for descriptor systems.
\newblock {\em {IEE} Proceedings - Control Theory and Applications},
  145(3):316--322, 1998.
\newblock \href {https://doi.org/10.1049/ip-cta:19982048}
  {\path{doi:10.1049/ip-cta:19982048}}.

\bibitem{Wei16}
H.~K. Weichelt.
\newblock {\em Numerical Aspects of Flow Stabilization by {R}iccati Feedback}.
\newblock {D}issertation, Department of Mathematics, Otto von Guericke
  University, Magdeburg, Germany, 2016.
\newblock \href {https://doi.org/10.25673/4493} {\path{doi:10.25673/4493}}.

\bibitem{Wil96}
C.~H.~K. Williamson.
\newblock Vortex dynamics in the cylinder wake.
\newblock {\em Annu. Rev. Fluid Mech.}, 28(1):477--539, 1996.
\newblock \href {https://doi.org/10.1146/annurev.fl.28.010196.002401}
  {\path{doi:10.1146/annurev.fl.28.010196.002401}}.

\bibitem{Zam81}
G.~Zames.
\newblock Feedback and optimal sensitivity: {M}odel reference transformations,
  multiplicative seminorms, and approximate inverses.
\newblock {\em {IEEE} Trans. Autom. Control}, 26(2):301--320, 1981.
\newblock \href {https://doi.org/10.1109/TAC.1981.1102603}
  {\path{doi:10.1109/TAC.1981.1102603}}.

\bibitem{ZhoDG96}
K.~Zhou, J.~C. Doyle, and K.~Glover.
\newblock {\em Robust and Optimal Control}.
\newblock Prentice-Hall, Upper Saddle River, NJ, 1996.

\end{thebibliography}

\end{document}